\documentclass{amsart}

\newcommand{\docclass}{amsart}

	\usepackage[T1]{fontenc}
\usepackage[latin9]{inputenc}

\usepackage{amssymb}
\usepackage{amstext}
\usepackage{amsmath}
\usepackage[all]{xy}

\usepackage{ifthen}
\usepackage{xstring}

		\IfStrEqCase{\docclass}{%
			{amsart}{\usepackage{amsthm}}%
			{elsart}{\usepackage{amsthm}}%
			{springer}{\smartqed}%
		}

\usepackage[numbers]{natbib}

\allowdisplaybreaks  		
	\newcommand{\defineamstheorems}{
	\theoremstyle{plain}
	\newtheorem{thm}{Theorem}[section]

	\theoremstyle{plain}
	\newtheorem{lem}[thm]{Lemma}

	\theoremstyle{definition}
	\newtheorem{defn}[thm]{Definition}

	\theoremstyle{remark}
	\newtheorem{rem}[thm]{Remark}

	\theoremstyle{plain}
	\newtheorem{prop}[thm]{Proposition}

	\theoremstyle{plain}
	\newtheorem{cor}[thm]{Corollary}

	\theoremstyle{definition}
	\newtheorem{exampleenv}[thm]{Example}

	\newenvironment{proofof}{\begin{proof}}{\end{proof}}
}

\newcommand{\definespringertheorems}{

	\newenvironment{thm}{\begin{theorem}}{\end{theorem}}
	\newenvironment{lem}{\begin{lemma}}{\end{lemma}}
	\newenvironment{defn}{\begin{definition}}{\end{definition}}
	\newenvironment{rem}{\begin{remark}}{\end{remark}}
	\newenvironment{prop}{\begin{proposition}}{\end{proposition}}
	\newenvironment{cor}{\begin{corollary}}{\end{corollary}}
	\newenvironment{exampleenv}{\begin{example}}{\end{example}}

	\newenvironment{proofof}{\begin{proof}}{\qed\end{proof}}

}

		\IfStrEqCase{\docclass}{%
			{amsart}{\defineamstheorems}%
			{elsart}{\defineamstheorems}%
			{springer}{\definespringertheorems}%
		}

\begin{document}
		
		\newcommand{\doctitle}{Normality of spaces of operators and quasi-lattices}
\newcommand{\shortdoctitle}{\doctitle}

\newcommand{\docabstract}{%
	We give an overview of normality and conormality properties of pre-ordered
	Banach spaces. For pre-ordered Banach spaces $X$ and $Y$ with closed
	cones we investigate normality of $B(X,Y)$ in terms of normality
	and conormality of the underlying spaces $X$ and $Y$. 

	Furthermore, we define a class of ordered Banach spaces called quasi-lattices
	which strictly contains the Banach lattices, and we prove that every
	strictly convex reflexive ordered Banach space with a closed proper
	generating cone is a quasi-lattice. These spaces provide a large class
	of examples $X$ and $Y$ that are not Banach lattices, but for which
	$B(X,Y)$ is normal. In particular, we show that a Hilbert space $\mathcal{H}$
	endowed with a Lorentz cone is a quasi-lattice (that is not a Banach
	lattice if $\dim\mathcal{H}\geq3$), and satisfies an identity analogous
	to the elementary Banach lattice identity $\||x|\|=\|x\|$ which holds
	for all elements $x$ of a Banach lattice. This is used to show that
	spaces of operators between such ordered Hilbert spaces are always
	absolutely monotone and that the operator norm is positively attained,
	as is also always the case for spaces of operators between Banach
	lattices.
}

		\IfStrEqCase{\docclass}{%
			{amsart}{\title{Normality of spaces of operators and quasi-lattices}}%
			{elsart}{\title{\doctitle}}%
			{springer}{\title{\doctitle}}%
		}
		\begin{abstract}%
			\docabstract%
		\end{abstract}


		\newcommand{\authornameMiek}{Miek Messerschmidt}
\newcommand{\emailMiek}{mmesserschmidt@gmail.com}

\newcommand{\addressLeiden}{%
	Mathematical Institute, 
	Leiden University, 
	P.O. Box 9512, 
	\mbox{2300 RA} Leiden, 
	The Netherlands
}

\newcommand{\addressNWU}{
	Unit for BMI,
	North-west University,
	Private Bag X6001,
	Potchefstroom,
	South Africa,
	2520
}

\newcommand{\authornameMarcel}{Marcel de Jeu}
\newcommand{\emailMarcel}{mdejeu@math.leidenuniv.nl}

\newcommand{\makeamsbio}{%
	\author{\authornameMiek}%
	\address{\authornameMiek, \addressLeiden}%
	\address{\footnote{Current address}\authornameMiek, \addressNWU}%
	\email{\emailMiek}
}

\newcommand{\makeelsbio}{%
	\address[leiden]{\addressLeiden}

	\author[leiden]{\authornameMarcel}
	\ead{\emailMarcel}
	
	\author[leiden]{\authornameMiek\corref{cor1}\fnref{fn1}}
	\ead{\emailMiek}
	\cortext[cor1]{Corresponding author}
	
	\fntext[fn1]{Current address: \addressNWU}
}

\newcommand{\makespringerbio}{%

	\author{\authornameMiek}

	\institute{%
		\authornameMiek %
			\at %
				\addressLeiden %
			\and 	%
				Current Address: \addressNWU %
			\\\email{\emailMiek} %
	}
}

		\newcommand{\subjectClassesForThisPaper}{%
		Primary: 46B40; %
		Secondary: 	 47B60\subjclasssep%
					 47H07\subjclasssep%
					 46B42\subjclasssep%
					 46A40\ignorespaces%
}

\newcommand{\keywordsForThisPaper}{
	(pre)-ordered Banach space\keywordssep%
	 operator norm\keywordssep%
	 quasi-lattice\keywordssep%
	 normality\keywordssep%
	 conormality\keywordssep%
	 Lorentz cone\ignorespaces%
}

\newcommand{\makeamsothermeta}{%
	\newcommand{\keywordssep}{, }
	\newcommand{\subjclasssep}{, }
	
	\keywords{\keywordsForThisPaper}
	\subjclass[2010]{\subjectClassesForThisPaper}
}

\newcommand{\makeelsothermeta}{%
	\newcommand{\keywordssep}{\sep }
	\newcommand{\subjclasssep}{\sep }

	\begin{keyword}%
		\keywordsForThisPaper%
		\MSC[2010]{\subjectClassesForThisPaper}%
	\end{keyword}%
}

\newcommand{\makespringerothermeta}{%
	\newcommand{\keywordssep}{\and }
	\newcommand{\subjclasssep}{\and }

	\keywords{\keywordsForThisPaper}

	\subclass{\subjectClassesForThisPaper}
}

		\IfStrEqCase{\docclass}{%
			{amsart}{\makeamsbio}%
			{elsart}{\makeelsbio}%
			{springer}{\makespringerbio}%
		}
	
		\IfStrEqCase{\docclass}{%
			{amsart}{\makeamsothermeta}%
			{elsart}{\makeelsothermeta}%
			{springer}{\makespringerothermeta}%
		}	
		\maketitle 
		
	
		\section{Introduction}

This paper's main aim is to investigate normality and monotonicity
(defined in Section \ref{sec:Normality-and-Conormality}) of pre-ordered
spaces of operators between pre-ordered Banach spaces. This investigation
is motivated by the relevance of this notion in the theory of positive
semigroups on pre-ordered Banach spaces \cite{BattyRobinson}, and
in the positive representation theory of groups and pre-ordered algebras
on pre-ordered Banach spaces \cite{CrossedprodsPartIII}. 

If $X$ and $Y$ are Banach lattices an elementary calculation shows
that the space $B(X,Y)$ is absolutely monotone, i.e., for $T,S\in B(X,Y)$,
if $\pm T\leq S$, then $\|T\|\leq\|S\|$. If $X$ and $Y$ are general
pre-ordered Banach spaces the situation is not so clear, and raises
a number of questions: If $B(X,Y)$ is, e.g., absolutely monotone,
does this necessarily imply that $X$ and $Y$ are Banach lattices?
If not, what are examples of pre-ordered Banach spaces $X$ and $Y$,
not being Banach lattices, such that $B(X,Y)$ is absolutely monotone?
What are the more general necessary and/or sufficient conditions $X$
and $Y$ have to satisfy for $B(X,Y)$ to be absolutely monotone?
This paper will attempt to answer such questions through an investigation
of the notions of normality and conormality of pre-ordered Banach
spaces which describe various ways in which cones interact with norms. 

A substantial part will devoted to introducing a class of ordered
Banach spaces, called quasi-lattices, which will furnish us with many
examples that are not necessarily Banach lattices. Quasi-lattices
occur in two slightly different forms, one of which includes all Banach
lattices (cf.\,Proposition \ref{prop:Banach-lattices-are-mu-quasi}).
We give a brief sketch of their construction. 

There are many pre-ordered Banach spaces with closed proper generating
cones that are not normed Riesz spaces, e.g., the finite dimensional
spaces $\mathbb{R}^{n}$ (with $n\geq3$) endowed with Lorentz cones
or endowed with polyhedral cones whose bases (in the sense of \cite[Section 1.7]{CnD})
are not $(n-1)$-simplexes. Although there is often an abundance of
upper bounds of arbitrary pairs of elements, none of them is a least
upper bound with respect to the ordering defined by the cone. An interesting
situation arises when one takes the norm into account when studying
the set of upper bounds of arbitrary pairs of elements. Even though
there might not exist a least upper bound with respect to the ordering
defined by the cone, there often exists a unique upper bound, called
the \emph{quasi-supremum}, which minimizes the sum of the distances
from this upper bound to the given two elements. This allows us to
define what will be called a \emph{quasi-lattice structure }on certain
ordered Banach spaces which might not be lattices (cf.\,Definition~\ref{def:quasi-lattices}). Surprisingly, many elementary vector lattice
properties for Riesz spaces carry over nearly verbatim to such spaces
(cf.\,Theorem \ref{thm:formulas1}), and in the case that a space
is a Banach lattice, its quasi-lattice structure and lattice structure
actually coincide (cf.\,Proposition \ref{prop:Banach-lattices-are-mu-quasi}).

Quasi-lattices occur in relative abundance, in fact, every strictly
convex reflexive ordered Banach space with a closed proper generating
cone is a quasi-lattice (cf.\,Theorem \ref{thm:reflexive-u-quasi-lattice}).
This will be used to show that every Hilbert space $\mathcal{H}$
endowed with a Lorentz cone is a quasi-lattice (which is not a Banach
lattice if $\dim(\mathcal{H})\geq3$). Such spaces will serve as examples
of spaces, which are not Banach lattices, such that the spaces of
operators between them are absolutely monotone (cf.\,Theorem \ref{thm:hilbert-lorentz-collected-properties}),
hence resolving the question of the existence of such spaces as posed
above.

We briefly describe the structure of the paper.

After giving preliminary definitions and terminology in Section \ref{sec:Preliminary-definitions},
we introduce various versions of the concepts of normality and conormality
of pre-ordered Banach spaces with closed cones in Section \ref{sec:Normality-and-Conormality}.
Normality is a more general notion than monotonicity, and roughly
is a measure of `the obtuseness/bluntness of a cone' (with respect
to the norm). Conormality roughly is a measure of `the acuity/sharpness
of a cone' (with respect to the norm). Normality and conormality properties
often occur in dual pairs, where a pre-ordered Banach space with a
closed cone has a normality property precisely when its dual has the
appropriate conormality property (cf.\,Theorem \ref{thm:normality-conormality-duality-properties}).
The terms `monotonicity' and `normality' are fairly standard throughout
the literature. However, the concept of conormality occurs scattered
under many names throughout the literature (chronologically, \cite{GrosbergKrein,Bonsall,Ando,Ellis,Davies,Ng,Wickstead,NgOpenDecompProperty,WongNg,Walsh,RobinsonYamamuro,Yamamuro,BattyRobinson,Clement,NgLaw}).
Although the definitions and results in Section \ref{sec:Normality-and-Conormality}
are not new, they are collected here in an attempt to give an overview
and to standardize the terminology.

In Section \ref{sec:Normality-of-operator-spaces}, with $X$ and
$Y$ pre-ordered Banach spaces with closed cones, we investigate the
normality of $B(X,Y)$ in terms of the normality and conormality of
$X$ and $Y$. Roughly, excluding degenerate cases, some form of conormality
of $X$ and normality of $Y$ is necessary and sufficient for having
some form of normality of the pre-ordered Banach space $B(X,Y)$ (cf.\,Theorems
\ref{thm-yamamuro-results} and  \ref{thm-(co)normality-conditions-onXandY-implies-normality-ir-B(X,Y)}).
Again, certain results are not new, but are included for the sake
of completeness.

In Section \ref{sec:Quasi-lattices-basic-properties} we introduce
quasi-lattices, a class of pre-ordered Banach spaces spaces that strictly
includes the Banach lattices. We establish their basic properties,
in particular, basic vector lattice identities which carry over from
Riesz spaces to quasi-lattices (cf.\,Theorem \ref{thm:formulas1}).

In Section \ref{sec:More-on-quasi-lattices} we prove one of our main
results: Every strictly convex reflexive pre-ordered Banach space
with a closed proper and generating cone is a quasi-lattice. Hence
there are many quasi-lattices.

Finally, in Section \ref{sec:Quasi-lattices-with-positively-attained-operator-norms},
we show that real Hilbert spaces endowed with Lorentz cones are quasi-lattices
and satisfy an identity analogous to the elementary identity $\||x|\|=\|x\|$
which holds for all elements $x$ of a Banach lattice. This is used
to show, for real Hilbert spaces $\mathcal{H}_{1}$ and $\mathcal{H}_{2}$
endowed with Lorentz cones, that $B(\mathcal{H}_{1},\mathcal{H}_{2})$
is absolutely monotone.

\section{Preliminary definitions and notation \label{sec:Preliminary-definitions}}

Let $X$ be a Banach space over the real numbers. Its topological
dual will be denoted by $X'$. A subset $C\subseteq X$ will be called
a \emph{cone} if $C+C\subseteq C$ and $\lambda C\subseteq C$ for
all $\lambda\geq0$. If a cone $C$ satisfies $C\cap(-C)=\{0\}$,
it will be called a \emph{proper} \emph{cone}, and if $X=C-C$, it
will be said to be \emph{generating }(\emph{in $X$})\emph{. }
\begin{defn}
A pair $(X,C)$, with $X$ a Banach space and $C\subseteq X$ a cone,
will be called a \emph{pre-ordered Banach space. }If $C$ is a proper
cone, $(X,C)$\emph{ }will be called an \emph{ordered Banach space}.
We will often suppress explicit mention of the pair and merely say
that $X$ is a (pre-)ordered Banach space. When doing so, we will
denote the implicit cone by $X_{+}$ and refer to it as \emph{the
cone of }$X$. For any $x,y\in X$, by $x\geq y$ we will mean $x-y\in X_{+}$.
We do not exclude the possibilities $X_{+}=\{0\}$ or $X_{+}=X$,
and we do not assume that $X_{+}$ is closed.
\end{defn}
Let $X$ and $Y$ be pre-ordered Banach spaces. The space of bounded
linear operators from $X$ to $Y$ will be denoted by $B(X,Y)$ and
by $B(X)$ if $X=Y$. Unless otherwise mentioned, $B(X,Y)$ is always
endowed with the operator norm. The space $B(X,Y)$ is easily seen
to be a pre-ordered Banach space when endowed with the cone\emph{
}$B(X,Y)_{+}:=\{T\in B(X,Y):TX_{+}\subseteq Y_{+}\}$. In particular,
the topological dual $X'$ also becomes a pre-ordered Banach space
when endowed with the \emph{dual cone} $X_{+}':=B(X,\mathbb{R})_{+}$.
For any $f\in X'$ and $y\in Y$, we will define the operator $f\otimes y\in B(X,Y)$
by $(f\otimes y)(x):=f(x)y$ for all $x\in X$. It is easily seen
that $\|f\otimes y\|=\|f\|\|y\|$.

\section{Normality and Conormality\label{sec:Normality-and-Conormality}}

In the current section we will define some of the possible norm-cone
interactions that may occur in pre-ordered Banach spaces, and investigate
how they relate to norm-cone interactions in the dual. Historically,
these properties have been assigned to either the norm or the cone
(e.g., `a cone is normal' and `a norm is monotone'). We will not follow
this convention and rather assign these labels to the pre-ordered
Banach space as a whole to emphasize the norm-cone interaction. 

We attempt to collect all known results and to standardize the terminology.
The definitions and results in the current section are essentially
known, but are scattered throughout the literature under quite varied
terminology%
\footnote{A note on terminology: The terms `normality' (due to Krein \cite{KreinNormality})
and `monotonicity' are fairly standard terms throughout the literature.
Our consistent use of the adjective `absolute' is inspired by \cite{WongNg}
and mimics its use in the term `absolute value'. 

The concept that we will call `conormality' has seen numerous equivalent
definitions and the nomenclature is rather varied in the existing
literature. The term `conormality' is due to Walsh \cite{Walsh},
who studied the property in the context of locally convex spaces.
What we will call `$1$-max-conormality' occurs under the name `strict
bounded decomposition property' in \cite{Bonsall_Linear_operators_in_complete_positive_cones}.
The properties that we will call `approximate $1$-absolute conormality'
and `approximate $1$-conormality', were first defined (but not named)
respectively by Davies \cite{Davies} and Ng \cite{Ng}. Batty and
Robinson give equivalent definitions for our conormality properties
which they call `dominating' and `generating' \cite{BattyRobinson}.%
}. References are provided when known to the author. 
\begin{defn}
Let $X$ be a pre-ordered Banach space with a closed cone and $\alpha>0$. 

We define the following \textbf{\emph{normality properties}}\emph{:}
\begin{enumerate}
\item We will say $X$ is \emph{$\alpha$-max-normal} if, for any $x,y,z\in X$,
$z\leq x\leq y$ implies $\|x\|\leq\alpha\max\{\|y\|,\|z\|\}$.
\item We will say $X$ is \emph{$\alpha$-sum-normal} if, for any $x,y,z\in X$,
$z\leq x\leq y$ implies $\|x\|\leq\alpha(\|y\|+\|z\|)$.
\item We will say $X$ is \emph{$\alpha$-absolutely normal} if, for any
$x,y\in X$, $\pm x\leq y$ implies $\|x\|\leq\alpha\|y\|$. We will
say $X$ is \emph{absolutely monotone} if it is $1$-absolutely normal.
\item We will say $X$ is \emph{$\alpha$-normal} if, for any $x,y\in X$,
$0\leq x\leq y$ implies $\|x\|\leq\alpha\|y\|$. We will say $X$
is \emph{monotone} if it is $1$-normal.
\end{enumerate}
We define the following\textbf{\emph{ conormality properties}}\emph{:}
\begin{enumerate}
\item We will say $X$ is \emph{$\alpha$-sum-conormal} if, for any $x\in X$,
there exist some $a,b\in X_{+}$ such that $x=a-b$ and $\|a\|+\|b\|\leq\alpha\|x\|$.
We will say $X$ is \emph{approximately $\alpha$-sum-conormal} if,
for any $x\in X$ and $\varepsilon>0$, there exist some $a,b\in X_{+}$
such that $x=a-b$ and $\|a\|+\|b\|<\alpha\|x\|+\varepsilon$. 
\item We will say $X$ is \emph{$\alpha$-max-conormal} if, for any $x\in X$,
there exist some $a,b\in X_{+}$ such that $x=a-b$ and $\max\{\|a\|,\|b\|\}\leq\alpha\|x\|$.
We will say $X$ is \emph{approximately $\alpha$-max-conormal} if,
for any $x\in X$ and $\varepsilon>0$, there exist some $a,b\in X_{+}$
such that $x=a-b$ and $\max\{\|a\|,\|b\|\}<\alpha\|x\|+\varepsilon$. 
\item We will say $X$ is \emph{$\alpha$-absolutely conormal} if, for any
$x\in X$, there exist some $a\in X_{+}$ such that $\pm x\leq a$
and $\|a\|\leq\alpha\|x\|$. We will say $X$ is \emph{approximately
$\alpha$-absolutely conormal} if, for any $x\in X$ and $\varepsilon>0$,
there exist some $a\in X_{+}$ such that $\pm x\leq a$ and $\|a\|<\alpha\|x\|+\varepsilon$. 
\item We will say $X$ is \emph{$\alpha$-conormal} if, for any $x\in X$,
there exist some $a\in X_{+}$ such that $0,x\leq a$ and $\|a\|\leq\alpha\|x\|$.
We will say $X$ is \emph{approximately $\alpha$-conormal} if, for
any $x\in X$ and $\varepsilon>0$, there exist some $a\in X_{+}$
such that $\{0,x\}\leq a$ and $\|a\|<\alpha\|x\|+\varepsilon$. 
\end{enumerate}
\end{defn}
The following two results show the relationship between different
(co)normality properties and for the most part are immediate from
the definitions.
\begin{prop}
\label{prop:normality-relationships}For any fixed $\alpha>0$, the
following implications hold between normality properties of a pre-ordered
Banach space $X$ with a closed cone:

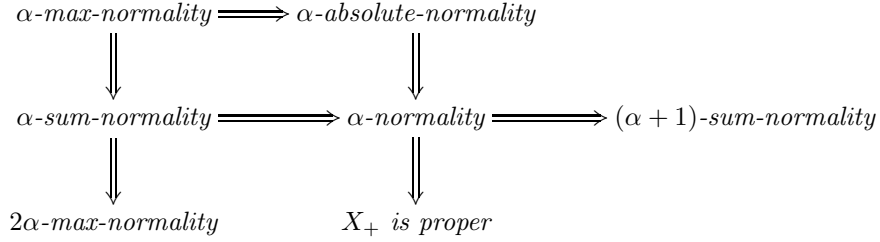
\[
\xymatrix{\alpha\mbox{-max-normality}\ar@2[d]\ar@2[r] & \alpha\mbox{-absolute-normality}\ar@2[d]\\
\alpha\mbox{-sum-normality}\ar@2[d]\ar@2[r] & \alpha\mbox{-normality}\ar@2[d]\ar@2[r] & (\alpha+1)\mbox{-sum-normality}\\
2\alpha\mbox{-max-normality} & X_{+}\mbox{ is proper}
}
\]
\end{prop}
\begin{proofof}
The only implication that is not immediate from the definitions is
that $\alpha$-normality implies $(\alpha+1)$-sum-normality. As to
this, let $X$ be an $\alpha$-normal pre-ordered Banach space with
a closed cone and $x,y,z\in X$ such that $z\leq x\leq y$. Then $0\leq x-z\leq y-z$,
so that, by $\alpha$-normality and the reverse triangle inequality,
\[
\|x\|-\|z\|\leq\|x-z\|\leq\alpha\|y-z\|\leq\alpha(\|y\|+\|z\|).
\]
Hence $\|x\|\leq\alpha\|y\|+(\alpha+1)\|z\|\leq(\alpha+1)(\|y\|+\|z\|)$.
\end{proofof}
Similar relationships hold between conormality properties as do between
normality properties. All implications follow immediately from the
definitions, with one exception. That is, if a pre-ordered Banach
space has a closed generating cone, then there exists a constant $\beta>0$
such that it is $\beta$-max-conormal (and hence $2\beta$-sum-conormal).
This is a result due to And\^o \cite[Lemma 1]{Ando}, and is the
bottom implication in the following proposition (although \cite[Lemma 1]{Ando}
assumes the cone to be proper, this is not necessary for its statement
to hold, cf.\,Theorem \ref{thm:continuous-upper-bound-function}).
\begin{prop}
\label{prop:conormality-relationships}For any fixed $\alpha>0$,
the following implications hold between conormality properties of
a pre-ordered Banach space $X$ with a closed cone:

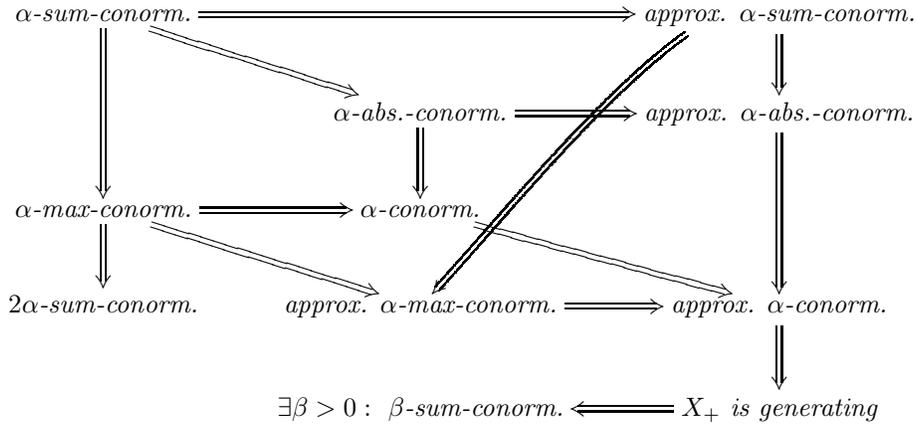
\[
\xymatrix{\alpha\mbox{-sum-conorm.}\ar@2[dd]\ar@2[rr]\ar@2[dr] &  & \mbox{approx. }\alpha\mbox{-sum-conorm.}\ar@2[d]\ar@2@(ul,ur)[lddd]\\
 & \alpha\mbox{-abs.-conorm.}\ar@2[d]\ar@2[r] & \mbox{approx. }\alpha\mbox{-abs.-conorm.}\ar@2[dd]\\
\alpha\mbox{-max-conorm.}\ar@2[d]\ar@2[dr]\ar@2[r] & \alpha\mbox{-conorm.}\ar@2[dr]\\
2\alpha\mbox{-sum-conorm.} & \mbox{approx. }\alpha\mbox{-max-conorm.}\ar@2[r] & \mbox{approx. }\alpha\mbox{-conorm.}\ar@2[d]\\
 & \exists\beta>0:\ \beta\mbox{-sum-conorm.} & X_{+}\mbox{ is generating}\ar@2[l]
}
\]

\end{prop}

\begin{rem}
The direct analogue to And\^o's Theorem \cite[Lemma 1]{Ando} (the
bottom implication in Proposition \ref{prop:conormality-relationships})
in Proposition \ref{prop:normality-relationships} would be that having
$X_{+}$ proper implies that $X$ is $\beta$-max-normal for some
$\beta>0$. This is false. Example \ref{exa:non-normal} gives a space
which has a proper cone but is not $\alpha$-normal for any $\alpha>0$. 
\end{rem}
{}
\begin{rem}
For the sake of completeness, we note that And\^o's Theorem \cite[Lemma 1]{Ando}
(the bottom implication in Proposition \ref{prop:conormality-relationships})
can be improved, in that the decomposition of elements into a difference
of elements from the cone can be chosen in a continuous, as well as
bounded and positively homogeneous manner. The following result is
a special case of \cite[Theorem 4.1]{RightInverses}, which is a general
principle for Banach spaces that are the sum of (not necessarily countably
many) closed cones. Its proof proceeds through an application of Michael's
Selection Theorem \cite[Theorem 17.66]{AliprantisBorder} and a generalization
of the usual Open Mapping Theorem \cite[Theorem 3.2]{RightInverses}:
\begin{thm}
\label{thm:continuous-upper-bound-function}Let $X$ be a pre-ordered
Banach space with a closed generating cone. Then there exist continuous
positively homogeneous functions $(\cdot)^{\pm}:X\to X_{+}$ and a
constant $\alpha>0$ such that $x=x^{+}-x^{-}$ and $\|x^{\pm}\|\leq\alpha\|x\|$
for all $x\in X$.
\end{thm}
\end{rem}
Normality and conormality properties often appear in dual pairs. Roughly,
a pre-ordered Banach space has a normality property if and only if
its dual has a corresponding conormality property, and vice versa.
The following theorem provides an overview of these normality-conormality
duality relationships as known to the author.
\begin{thm}
\label{thm:normality-conormality-duality-properties}Let $X$ be a
pre-ordered Banach space with a closed cone.
\begin{enumerate}
\item The following equivalences hold:

\begin{enumerate}
\item For $\alpha>0$, the space $X$ is $\alpha$-max-normal if and only
if $X'$ is $\alpha$-sum-conormal.
\item For $\alpha>0$, the space $X$ is $\alpha$-sum-normal if and only
if $X'$ is $\alpha$-max-conormal.
\item For $\alpha>0$, the space $X$ is $\alpha$-absolutely normal if
and only if $X'$ is $\alpha$-absolutely conormal.
\item For $\alpha>0$, the space $X$ is $\alpha$-normal if and only if
$X'$ is $\alpha$-conormal.
\item There exists an $\alpha>0$ such that $X$ is $\alpha$-max-normal
if and only if $X_{+}'$ is generating.
\end{enumerate}
\item The following equivalences hold:

\begin{enumerate}
\item For $\alpha>0$, the space $X$ is approximately $\alpha$-sum-conormal
if and only if $X'$ is $\alpha$-max-normal.
\item For $\alpha>0$, the space $X$ is approximately $\alpha$-max-conormal
if and only if $X'$ is $\alpha$-sum-normal.
\item For $\alpha>0$, the space $X$ is approximately $\alpha$-absolutely
conormal if and only if $X'$ is $\alpha$-absolutely normal.
\item For $\alpha>0$, the space $X$ is approximately $\alpha$-conormal
if and only if $X'$ is $\alpha$-normal.
\item The cone $X_{+}$ is generating if and only if there exists an $\alpha>0$
such that $X'$ is $\alpha$-max-normal.
\end{enumerate}
\end{enumerate}
\end{thm}
The result (1)(a) was first proven by Grosberg and Krein in \cite{GrosbergKrein}
(via \cite[Theorem 7]{Ellis}). The result (2)(a) was established
by Ellis \cite[Theorem 8]{Ellis}. For $\alpha=1$, the results (1)(c),(d),
(2)(c) and (d) are due to Ng \cite[Propositions 5, 6; Theorems 6, 7]{Ng}.
The fully general results (1)(d) and (2)(d) appear first in \cite[Theorem 1.1]{RobinsonYamamuro}
by Robinson and Yamamuro, and later in \cite[Theorems 1,2]{NgLaw}
by Ng an Law. Proofs of (1)(a) (again), (1)(b), (1)(c), (1)(d) (again),
and (2)(a) (again), (2)(b), (2)(c), and (2)(d) (again) are due to
Batty and Robinson in \cite[Theorems 1.1.4, 1.3.1, 1.2.2]{BattyRobinson}.
The results (1)(e) and (2)(e) are due to And\^o \cite[Theorem 1]{Ando}.

Bonsall proved an analogous duality result for locally convex spaces
in \cite[Theorem 2]{Bonsall}.

The following lemma shows that conormality properties and approximate
conormality properties of dual spaces are equivalent. Ng proved (3)
for the case $\alpha=1$ in \cite[Theorem 6]{Ng}:
\begin{lem}
\label{lem:dual-approx-conorm-equivalent-to-conorm}Let $X$ be a
pre-ordered Banach space with a closed cone. Then the following equivalences
hold:
\begin{enumerate}
\item For $\alpha>0$, the space  $X'$ is approximately $\alpha$-sum-conormal
if and only if $X'$ is $\alpha$-sum-conormal.
\item For $\alpha>0$, the space $X'$ is approximately $\alpha$-max-conormal
if and only if $X'$ is $\alpha$-max-conormal.
\item For $\alpha>0$, the space $X'$ is approximately $\alpha$-absolutely
conormal if and only if $X'$ is $\alpha$-absolutely conormal.
\item For $\alpha>0$, the space $X'$ is approximately $\alpha$-conormal
if and only if $X'$ is $\alpha$-conormal.
\end{enumerate}
\end{lem}
\begin{proofof}
That a conormality property implies the associated approximate conormality
property is trivial. We will therefore only prove the forward implications.

We prove (1). Let $X'$ be approximately $\alpha$-sum-conormal. Then,
for any $\beta>\alpha$ and any $0\neq f\in X$, by taking $\varepsilon=(\beta-\alpha)\|f\|>0$,
we have that there exist $g,h\in X'_{+}$ such that $f=g-h$ and $\|g\|+\|h\|\leq\alpha\|f\|+(\beta-\alpha)\|f\|=\beta\|f\|$.
Therefore, $X'$ is $\beta$-sum-conormal for every $\beta>\alpha$.
Now, by part (1)(a) of Theorem \ref{thm:normality-conormality-duality-properties},
$X$ is $\beta$-max-normal for every $\beta>\alpha$. Therefore,
if $x,y,z\in X$ are such that $z\leq x\leq y$, then $\|x\|\leq\beta\max\{\|y\|,\|z\|\}$
for all $\beta>\alpha$, and hence $\|x\|\leq\inf_{\beta>\alpha}\beta\max\{\|y\|,\|z\|\}=\alpha\max\{\|y\|,\|z\|\}$.
We conclude that $X$ is $\alpha$-max-normal, and, again by part
(1)(a) Theorem \ref{thm:normality-conormality-duality-properties},
that $X'$ is $\alpha$-sum-conormal.

The assertions (2), (3) and (4) follow through similar arguments.
\end{proofof}
By Theorem \ref{thm:normality-conormality-duality-properties} and
Lemma \ref{lem:dual-approx-conorm-equivalent-to-conorm}, a pre-ordered
Banach space with a closed cone possesses both a normality property
and its paired approximate conormality property (with the same constant)
if and only if its dual possesses the same properties (cf.\,Corollary
\ref{cor:regularity-duality}). Such spaces are called regular and
were first studied by Davies in \cite{Davies} and Ng in \cite{Ng}.
\begin{defn}
\label{def:regularity}Let $X$ be a pre-ordered Banach space with
a closed cone. We define the following \textbf{\emph{regularity properties}}\emph{:}%
\footnote{The term `regularity' is due to Davies \cite{Davies}. Our naming
convention is to attach the names of the persons who (to the author's
knowledge) first proved the relevant normality-conormality duality
results of the defining properties (cf.\,Theorem \ref{thm:normality-conormality-duality-properties}).%
}
\begin{enumerate}
\item For $\alpha>0$, we will say $X$ is \emph{$\alpha$-Ellis-Grosberg-Krein
regular} if $X$ is both $\alpha$-max-normal and approximately $\alpha$-sum-conormal\emph{.}
\item For $\alpha>0$, we will say $X$ is \emph{$\alpha$-Batty-Robinson
regular} if $X$ is both $\alpha$-sum-normal and approximately $\alpha$-max-conormal\emph{.}
\item For $\alpha>0$, we will say $X$ is \emph{$\alpha$-absolutely Davies-Ng
regular} if $X$ is both $\alpha$-absolutely normal and approximately
$\alpha$-absolutely conormal\emph{.}
\item For $\alpha>0$, we will say $X$ is \emph{$\alpha$-Davies-Ng regular}
if $X$ is both $\alpha$-normal and approximately $\alpha$-conormal\emph{.}
\item We will say $X$ is \emph{And\^o regular} if\emph{ $X_{+}$ }is generating
and there exists an $\alpha>0$ such that $X$ is $\alpha$-max-normal\emph{.}
\end{enumerate}
\end{defn}
It should be noted that every Banach lattice is $1$-absolutely Davies-Ng
regular. 

The following result combines Propositions \ref{prop:normality-relationships}
and \ref{prop:conormality-relationships} to provide relationships
that exist between regularity properties.
\begin{prop}
\label{prop:regularity-relationships}For any fixed $\alpha>0$, the
following implications hold between regularity properties of a pre-ordered
Banach space with a closed cone:

\[
\xymatrix{\alpha\mbox{-Ellis-Grosberg-Krein regularity}\ar@2[d]\ar@2[r] & \alpha\mbox{-Batty-Robinson regularity}\ar@2[d]\\
\alpha\mbox{-absolute Davies-Ng regularity}\ar@2[r] & \alpha\mbox{-Davies-Ng regularity}\ar@2[d]\\
\exists\beta>0:\ \beta\mbox{-Ellis-Grosberg-Krein regularity} & \mbox{And\^o
regularity}\ar@2[l] }
\]
\end{prop}
\begin{proofof}
The only implication that does not follow immediately from Propositions
\ref{prop:normality-relationships} and \ref{prop:conormality-relationships},
is that And\^o regularity implies $\beta$-Ellis-Grosberg-Krein regularity
for some $\beta>0$. As to this, let $X$ be an And\^o regular ordered
Banach space with a closed cone. By Proposition \ref{prop:conormality-relationships},
since $X_{+}$ is generating, there exists some $\delta>0$ such that
$X$ is $\delta$-sum-conormal. By assumption, there exists an $\alpha>0$,
such that $X$ is $\alpha$-max-normal. By taking $\beta:=\max\{\delta,\alpha\}$,
we see that $X$ is also $\beta$-max-normal and (approximately) $\beta$-sum-conormal.
We conclude that $X$ is $\beta$-Ellis-Grosberg-Krein regular.
\end{proofof}
A straightforward application of Theorem \ref{thm:normality-conormality-duality-properties}
and Lemma \ref{lem:dual-approx-conorm-equivalent-to-conorm} then
yields:
\begin{cor}
\label{cor:regularity-duality}Let $X$ be a pre-ordered Banach space
with a closed cone. Then the following equivalences hold:
\begin{enumerate}
\item For $\alpha>0$, the space $X$ is $\alpha$-Ellis-Grosberg-Krein
regular if and only if $X'$ is $\alpha$-Ellis-Grosberg-Krein regular.
\item For $\alpha>0$, the space $X$ is $\alpha$-\textup{\emph{Batty-Robinson}}
regular if and only if $X'$ is $\alpha$-\textup{\emph{Batty-Robinson
}}regular.
\item For $\alpha>0$, the space $X$ is $\alpha$-\textup{\emph{absolutely
Davies-Ng}}\emph{ }regular if and only if $X'$ is $\alpha$-\textup{\emph{absolutely
Davies-Ng}} regular.
\item For $\alpha>0$, the space $X$ is $\alpha$-\textup{\emph{Davies-Ng}}
regular if and only if $X'$ is $\alpha$-\textup{\emph{Davies-Ng}}
regular.
\item The space $X$ is And\^o regular if and only if $X'$ is And\^o
regular.
\end{enumerate}
\end{cor}

\section{The normality of pre-ordered Banach spaces of bounded linear
operators\label{sec:Normality-of-operator-spaces}}

If $X$ and $Y$ are pre-ordered Banach spaces with closed cones,
we investigate necessary and sufficient conditions for the pre-ordered
Banach space $B(X,Y)$ to have a normality property. Where results
are known to the author from the literature, references are provided. 

We begin, in the following result, by investigating necessary conditions
for $B(X,Y)$ to have a normality property. Parts (2) and (3) in the
special case $X=Y$ and $\alpha=1$ in the following theorem are due
Yamamuro \cite[1.2--3]{Yamamuro}. Batty and Robinson also proved
part (2) for $X=Y$ and $\alpha=1$, and part (3) for $\alpha=\beta=1$
\cite[Corollary 1.7.5, Proposition 1.7.6]{BattyRobinson}. Part (5)
is due to Wickstead \cite[Theorem 3.1]{Wickstead}.
\begin{thm}
\label{thm-yamamuro-results}Let $X$ and $Y$ be non-zero pre-ordered
Banach spaces with closed cones and $\alpha>0$. 
\begin{enumerate}
\item The cone $B(X,Y)_{+}$ is proper if and only if $X=\overline{X_{+}-X_{+}}$
and $Y_{+}$ is proper.
\item Let $B(X,Y)$ be $\alpha$-normal. If $Y_{+}\neq\{0\}$, then $X$
is approximately $\alpha$-conormal. If $X_{+}'\neq\{0\}$, then $Y$
is $\alpha$-normal.
\item Let $B(X,Y)$ be $\alpha$-absolutely normal. If $Y_{+}\neq\{0\}$,
then $X$ is approximately $\alpha$-absolutely conormal. If $X_{+}'\neq\{0\}$,
then $Y$ is $\alpha$-absolutely normal.
\item Let $B(X,Y)$ be $\alpha$-sum-normal. If $Y_{+}\neq\{0\}$, then
$X$ is approximately $\alpha$-max-conormal. If $X_{+}'\neq\{0\}$,
then $Y$ is $\alpha$-sum-normal.
\item Let $B(X,Y)$ be $\alpha$-max-normal. If $Y_{+}\neq\{0\}$, then
$X$ is approximately $\alpha$-sum-conormal. If $X_{+}'\neq\{0\}$,
then $Y$ is $\alpha$-max-normal.
\end{enumerate}
\end{thm}
\begin{proofof}
We prove (1). Let $B(X,Y)_{+}$ be proper. Suppose $X\neq\overline{X_{+}-X_{+}}$.
By the Hahn-Banach Theorem there exists a non-zero functional $f\in X'$
such that $f|_{\overline{X_{+}-X_{+}}}=0$. Let $0\neq y\in Y$, then
$\pm f\otimes y\geq0$ since $f\otimes y|_{X_{+}}=0$. Therefore $B(X,Y)_{+}$
is not proper, contradicting our assumption. Suppose $Y_{+}$ is not
proper. Let $0\neq y\in Y_{+}\cap(-Y_{+})$ and $0\neq f\in X'$.
Then $\pm f\otimes y\geq0$, and hence $B(X,Y)_{+}$ is not proper,
contradicting our assumption.

Let $X=\overline{X_{+}-X_{+}}$ and $Y_{+}$ be proper. If $T\in B(X,Y)_{+}\cap(-B(X,Y)_{+})$,
then, since $Y_{+}$ is proper, $TX_{+}=\{0\}$. Hence $T(X_{+}-X_{+})=\{0\}$,
and by density of $X_{+}-X_{+}$ in $X$, we have $T=0$.

We prove (2). Let $B(X,Y)$ be $\alpha$-normal. With $Y_{+}\neq\{0\}$,
by Theorem \ref{thm:normality-conormality-duality-properties}, to
conclude that $X$ is approximately $\alpha$-conormal, it is sufficient
to prove that $X'$ is $\alpha$-normal. Let $f,g\in X'$ satisfy
$0\leq f\leq g$, and let $0\neq y\in Y_{+}$. Then $0\leq f\otimes y\leq g\otimes y$,
and by the $\alpha$-normality of $B(X,Y)$, 
\[
\|f\|\|y\|=\|f\otimes y\|\leq\alpha\|g\otimes y\|=\alpha\|g\|\|y\|.
\]
Therefore $\|f\|\leq\alpha\|g\|$, and hence $X'$ is $\alpha$-conormal.
With $X_{+}'\neq\{0\}$, let $0\neq f\in X_{+}'$ be arbitrary, and
$y,z\in Y$ such that $0\leq y\leq z$. Then $0\leq f\otimes y\leq f\otimes z$
in $B(X,Y)$, and by the $\alpha$-normality of $B(X,Y)$, 
\[
\|f\|\|y\|=\|f\otimes y\|\leq\alpha\|f\otimes z\|=\alpha\|f\|\|z\|.
\]
Hence $\|y\|\leq\alpha\|z\|$ and we conclude that $Y$ is $\alpha$-normal.

We prove (3). Let $B(X,Y)$ be $\alpha$-absolutely normal. With $Y_{+}\neq\{0\}$,
by Theorem \ref{thm:normality-conormality-duality-properties}, to
conclude that $X$ is approximately $\alpha$-absolutely conormal,
it is sufficient to prove that $X'$ is $\alpha$-absolutely normal.
Let $f,g\in X'$ satisfy $\pm f\leq g$, and let $0\neq y\in Y_{+}$.
Then $\pm f\otimes y\leq g\otimes y$, and by the $\alpha$-absolute
normality of $B(X,Y)$, 
\[
\|f\|\|y\|=\|f\otimes y\|\leq\alpha\|g\otimes y\|=\alpha\|g\|\|y\|.
\]
Therefore $\|f\|\leq\alpha\|g\|$, and hence $X'$ is $\alpha$-absolutely
normal. With $X_{+}'\neq\{0\}$, let $0\neq f\in X_{+}'$ be arbitrary,
and $y,z\in Y$ such that $\pm y\leq z$. Then $\pm f\otimes y\leq f\otimes z$
in $B(X,Y)$, and by the $\alpha$-absolute normality of $B(X,Y)$,
\[
\|f\|\|y\|=\|f\otimes y\|\leq\alpha\|f\otimes z\|=\alpha\|f\|\|z\|
\]
Hence $\|y\|\leq\alpha\|z\|$ and we conclude that $Y$ is $\alpha$-absolutely
normal.

We prove (4). Let $B(X,Y)$ is $\alpha$-sum-normal. With $Y_{+}\neq\{0\}$,
by Theorem \ref{thm:normality-conormality-duality-properties}, it
is sufficient to prove that $X'$ is $\alpha$-sum-normal to conclude
that $X$ is approximately $\alpha$-max-conormal. Let $0\neq y\in Y_{+}$
and $f,g,h\in X'$ satisfy $g\leq f\leq h$. Then $g\otimes y\leq f\otimes y\leq h\otimes y$
in $B(X,Y)$, and by the $\alpha$-sum-normality of $B(X,Y)$, 
\[
\|f\|\|y\|=\|f\otimes y\|\leq\alpha\left(\|g\otimes y\|+\|h\otimes y\|\right)=\alpha(\|g\|+\|h\|)\|y\|.
\]
Hence $\|f\|\leq\alpha(\|g\|+\|h\|)$ and $X'$ is $\alpha$-sum-normal.
With $X_{+}'\neq\{0\}$, to prove that $Y$ is $\alpha$-sum-normal,
let $u,v,y\in Y$ satisfy $u\leq y\leq v$ and let $0\neq f\in X_{+}'$.
Then $f\otimes u\leq f\otimes y\leq f\otimes v$ in $B(X,Y)$, and
hence, $\|y\|\leq\alpha(\|u\|+\|v\|)$ as before.

The proof of (5) is analogous to that of (4).
\end{proofof}
Converse-like implications to the previous result also hold, giving
sufficient conditions for $B(X,Y)$ to have a normality property.
Part (1) and the case $\alpha=\beta=1$ of part (3) are due to Batty
and Robinson \cite[Proposition 1.7.3, Corollary 1.7.5]{BattyRobinson}.
The special case $X=Y$ and $\alpha=\beta=1$ of part (3) is due to
Yamamuro \cite[1.3]{Yamamuro}. The case where $X$ is approximately
$\alpha$-sum-conormal and $Y$ is $\beta$-max-normal of part (4)
is due to Wickstead \cite[Theorem 3.1]{Wickstead}. 
\begin{thm}
\label{thm-(co)normality-conditions-onXandY-implies-normality-ir-B(X,Y)}Let
$X$ and $Y$ be pre-ordered Banach spaces with closed cones and $\alpha,\beta>0$. 
\begin{enumerate}
\item If $X_{+}$ is generating and $Y$ is $\alpha$-normal, then there
exists some $\gamma>0$ for which $B(X,Y)$ is $\gamma$-normal.
\item If $X$ is approximately $\alpha$-conormal and $Y$ is $\beta$-normal,
then $B(X,Y)$ is $(2\alpha+1)\beta$-normal.
\item If $X$ is approximately $\alpha$-absolutely conormal and $Y$ is
$\beta$-absolutely normal, then $B(X,Y)$ is $\alpha\beta$-absolutely
normal.
\item If $X$ is approximately $\alpha$-sum-conormal and $Y$ is $\beta$-normal
\textup{(}$\beta$-absolutely normal, $\beta$-max-normal, $\beta$-sum-normal
respectively\textup{)}, then $B(X,Y)$ is $\alpha\beta$-normal \textup{(}$\alpha\beta$-absolutely
normal, $\alpha\beta$-max-normal, $\alpha\beta$-sum-normal respectively\textup{)}
\end{enumerate}
\end{thm}
\begin{proofof}
We prove (1). By And\^o's Theorem \cite[Lemma 1]{Ando}, the fact
that $X_{+}$ is generating in $X$ implies that there exists some
$\beta>0$ such that $X$ is $\beta$-max-conormal. Let $T,S\in B(X,Y)$
be such that $0\leq T\leq S$. Then, for any $x\in X$, let $a,b\in X_{+}$
be such that $x=a-b$ and $\max\{\|a\|,\|b\|\}\leq\beta\|x\|$, so
that $0\leq Ta\leq Sa$ and $0\leq Tb\leq Sb$. By $\alpha$-normality
of $Y$, 
\[
\|Tx\|\leq\|Ta\|+\|Tb\|\leq\alpha(\|Sa\|+\|Sb\|)\leq\alpha\|S\|(\|a\|+\|b\|)\leq2\alpha\beta\|S\|\|x\|,
\]
hence $\|T\|\leq2\alpha\beta\|S\|$.

We prove (2). Let $T,S\in B(X,Y)$ be such that $0\leq T\leq S$.
Let $x\in X$ be arbitrary. Then, for every $\varepsilon>0$, there
exists some $a\in X_{+}$ such that $\{0,x\}\leq a$ and $\|a\|\leq\alpha\|x\|+\varepsilon$.
Since $x=a-(a-x)$ and $a,a-x\geq0$, we obtain $0\leq Ta\leq Sa$
and $0\leq T(a-x)\leq S(a-x)$, and hence
\begin{eqnarray*}
\|Tx\| & = & \|Ta-T(a-x)\|\\
 & \leq & \|Ta\|+\|T(a-x)\|\\
 & \leq & \beta\|Sa\|+\beta\|S(a-x)\|\\
 & \leq & \beta\|S\|(\alpha\|x\|+\varepsilon)+\beta\|S\|(\alpha\|x\|+\varepsilon+\|x\|)\\
 & = & (2\alpha+1)\beta\|S\|\|x\|+2\varepsilon\beta\|S\|.
\end{eqnarray*}
Since $\varepsilon>0$ was chosen arbitrarily, we conclude that $\|T\|\leq(2\alpha+1)\beta\|S\|$.

We prove (3). Let $T,S\in B(X,Y)$ satisfy $\pm T\leq S$. Let $x\in X$
be arbitrary. Then, for every $\varepsilon>0$, there exists some
$a\in X_{+}$ satisfying $\pm x\leq a$ and $\|a\|<\alpha\|x\|+\varepsilon$.
Then 
\[
Tx=T\left(\frac{a+x}{2}\right)-T\left(\frac{a-x}{2}\right),
\]
and hence, 
\[
\pm Tx=\pm T\left(\frac{a+x}{2}\right)\mp T\left(\frac{a-x}{2}\right).
\]
Since $(a+x)/2\geq0$, $(a-x)/2\geq0$ and $\pm T\leq S$, we find
\[
\pm Tx\leq S\left(\frac{a+x}{2}\right)+S\left(\frac{a-x}{2}\right)=Sa.
\]
Now, because $Y$ is $\beta$-absolutely normal, we obtain 
\[
\|Tx\|\leq\beta\|Sa\|\leq\beta\|S\|\|a\|\leq\alpha\beta\|S\|\|x\|+\varepsilon\beta\|S\|.
\]
Since $\varepsilon>0$ was chosen arbitrarily, we conclude that $B(X,Y)$
is $\alpha\beta$-absolutely normal.

We prove (4). Let $X$ be approximately $\alpha$-sum-conormal and
let $Y$ be $\beta$-normal. Let $T,U\in B(X,Y)$ satisfy $0\leq T\leq U$
and let $x\in X$ be arbitrary. Then, for every $\varepsilon>0$,
there exist $x_{1},x_{2}\in X_{+}$ such that $x=x_{1}-x_{2}$ and
$\|x_{1}\|+\|x_{2}\|<\alpha\|x\|+\varepsilon$. Also, $0\leq Tx_{i}\leq Ux_{i}$
implies $\|Tx_{i}\|\leq\beta\|Ux_{i}\|$ for $i=1,2$. Therefore,
\begin{eqnarray*}
\|Tx\| & \leq & \|Tx_{1}\|+\|Tx_{2}\|\\
 & \leq & \beta\|Ux_{1}\|+\beta\|Ux_{2}\|\\
 & \leq & \beta\|U\|(\|x_{1}\|+\|x_{2}\|)\\
 & \leq & \alpha\beta\|U\|\|x\|+\varepsilon\beta\|U\|.
\end{eqnarray*}
Since $x\in X$ and $\varepsilon>0$ were arbitrary, we may conclude
that $B(X,Y)$ is $\alpha\beta$-normal. The case where $X$ is approximately
$\alpha$-sum-conormal and $Y$ is $\beta$-absolutely normal follows
similarly.

Let $X$ be approximately $\alpha$-sum-conormal and let $Y$ be $\beta$-max-normal.
Let $T,U,V\in B(X,Y)$ satisfy $U\leq T\leq V$ and let $x\in X$
be arbitrary. Then, for every $\varepsilon>0$, there exist $x_{1},x_{2}\in X_{+}$
such that $x=x_{1}-x_{2}$ and $\|x_{1}\|+\|x_{2}\|<\alpha\|x\|+\varepsilon$.
Also, $Ux_{i}\leq Tx_{i}\leq Vx_{i}$ implies $\|Tx_{i}\|\leq\beta\max\{\|Ux_{i}\|,\|Vx_{i}\|\}$
for $i=1,2$. Therefore, 
\begin{eqnarray*}
\|Tx\| & \leq & \|Tx_{1}\|+\|Tx_{2}\|\\
 & \leq & \beta\max\{\|Ux_{1}\|,\|Vx_{1}\|\}+\beta\max\{\|Ux_{2}\|,\|Vx_{2}\|\}\\
 & \leq & \beta\max\{\|U\|,\|V\|\}(\|x_{1}\|+\|x_{2}\|)\\
 & \leq & \alpha\beta\max\{\|U\|,\|V\|\}\|x\|+\varepsilon\beta\max\{\|U\|,\|V\|\}.
\end{eqnarray*}
Since $x\in X$ and $\varepsilon>0$ were arbitrary, we may conclude
that $B(X,Y)$ is $\alpha\beta$-max-normal. The case where $X$ is
approximately $\alpha$-sum-conormal and $Y$ is $\beta$-sum-normal
follows similarly.
\end{proofof}
If one has further knowledge of the behavior of the positive bounded
linear operators, specifically that their norms are determined by
their behavior on the cone, then one can improve the constant in (2)
of the above theorem. This will be discussed in the rest of this section.
\begin{defn}
Let $X$ be a pre-ordered Banach space with a closed cone and $Y$
a Banach space. For $T\in B(X,Y)$, we define $\|T\|_{+}:=\sup\{\|Tx\|:x\in X_{+},\ \|x\|=1\}$. 

If $X=\overline{X_{+}-X_{+}},$ then $\|\cdot\|_{+}$ is a norm on
$B(X,Y)$, called the \emph{Robinson norm} (as named by Yamamuro in
\cite{Yamamuro}). We will say that the operator norm on $B(X,Y)$
is \emph{positively attained} (as named by Batty and Robinson in \cite{BattyRobinson})
if $\|T\|=\|T\|_{+}$ for all positive operators $T\in B(X,Y)_{+}$. 
\end{defn}
If $X_{+}$ is closed and generating, $\|\cdot\|_{+}$ is in fact
equivalent to the usual operator norm on $B(X,Y)$. The following
result is a slight refinement of a remark by Batty and Robinson \cite[p.\,248]{BattyRobinson}.
\begin{prop}
If $X$ is a pre-ordered Banach space with a closed generating cone
and $Y$ a Banach space, then the Robinson norm is equivalent to the
operator norm on $B(X,Y)$.\end{prop}
\begin{proofof}
By And\^o's Theorem \cite[Lemma 1]{Ando}, $X$ is $\alpha$-max-conormal
for some $\alpha>0$. Let $x\in X$ and $T\in B(X,Y)$ be arbitrary,
then there exist $a,b\in X_{+}$ such that $x=a-b$ and $\max\{\|a\|,\|b\|\}\leq\alpha\|x\|$.
Hence 
\begin{eqnarray*}
\|Tx\| & = & \|Ta-Tb\|\\
 & \leq & \|Ta\|+\|Tb\|\\
 & \leq & \|T\|_{+}(\|a\|+\|b\|)\\
 & \leq & 2\alpha\|T\|_{+}\|x\|.
\end{eqnarray*}
Therefore, $\|T\|_{+}\leq\|T\|\leq2\alpha\|T\|_{+}$.
\end{proofof}
Part (2) of Theorem \ref{thm-(co)normality-conditions-onXandY-implies-normality-ir-B(X,Y)}
can be improved if we know that the operator norm is positively attained.
\begin{prop}
\label{thm-pos-attained-opnorm-implies-normality}Let $X$ and $Y$
be pre-ordered Banach spaces with closed cones, with $Y$ $\alpha$-normal
for some $\alpha>0$. If the operator norm on $B(X,Y)$ is positively
attained, then $B(X,Y)$ is $\alpha$-normal. \end{prop}
\begin{proofof}
Let $T,S\in B(X,Y)$ satisfy $0\leq T\leq S$. Then, for any $x\in X_{+}$,
$0\leq Tx\leq Sx$, and hence $\|Tx\|\leq\alpha\|Sx\|$. We then see
that 
\begin{eqnarray*}
\|T\| & = & \|T\|_{+}\\
 & = & \sup\{\|Tx\|:x\in X_{+},\ \|x\|\leq1\}\\
 & \leq & \alpha\sup\{\|Sx\|:x\in X_{+},\ \|x\|\leq1\}\\
 & = & \alpha\|S\|_{+}\\
 & = & \alpha\|S\|,
\end{eqnarray*}
and conclude that $B(X,Y)$ is $\alpha$-normal.
\end{proofof}
The following theorem gives one necessary condition and some sufficient
conditions to have that an operator norm is positively attained. The
sufficiency of (1)%
\footnote{There is a small error in the statement of (1) in \cite[Proposition 1.7.8.]{BattyRobinson}.
We give its correct statement and proof.%
}, (2), and the necessity of approximate $1$-conormality in the following
theorem are due to Batty and Robinson in \cite[Proposition 1.7.8.]{BattyRobinson}. 
\begin{thm}
\label{thm:sufficient-conditions-for-positively-attained-operator-norm}Let
$X$ and $Y$ be pre-ordered Banach spaces with closed cones. 

If $Y_{+}\neq\{0\}$ and the operator norm on $B(X,Y)$ is positively
attained, then $X$ is approximately $1$-conormal.

Any of the following conditions is sufficient for the operator norm
on $B(X,Y)$ to be positively attained:
\begin{enumerate}
\item The space $X$ is approximately $1$-max-conormal and $Y$ is $1$-max-normal.
\item The space $X$ is approximately $1$-absolutely conormal and $Y$
is absolutely monotone \textup{(}i.e., if $X=Y$, $X$ is $1$-absolutely
Davies-Ng regular\textup{)}.
\item The space $X$ is approximately $1$-sum-conormal \textup{(}in which
case $\|T\|=\|T\|_{+}$ even holds for all $T\in B(X,Y)$\textup{)}.
\end{enumerate}
\end{thm}
\begin{proofof}
We prove the necessity of approximate $1$-conormality of $X$ when
$Y_{+}\neq\{0\}$ and the operator norm on $B(X,Y)$ is positively
attained. Let $f\in X_{+}'$ be arbitrary and let $0\neq y\in Y_{+}$.
Then, since the operator norm on $B(X,Y)$ is positively attained,
\[
\|f\|\|y\|=\|f\otimes y\|=\|f\otimes y\|_{+}=\|f\|_{+}\|y\|,
\]
so that $\|f\|=\|f\|_{+}$. For all $f,g\in X'$ satisfying $0\leq f\leq g$,
we obtain $\|f\|=\|f\|_{+}\leq\|g\|_{+}=\|g\|$. Therefore $X'$ is
monotone, and by part (2)(d) of Theorem \ref{thm:normality-conormality-duality-properties},
$X$ is approximately $1$-conormal.

We prove the sufficiency of (1). Let $T\in B(X,Y)_{+}$. Let $x\in X$
and $\varepsilon>0$ be arbitrary. Then, since $X$ is $1$-max-conormal,
there exist $a,b\in X_{+}$ such that $x=a-b$ and $\max\{\|a\|,\|b\|\}<\|x\|+\varepsilon$.
We notice that $-b\leq x\leq a$ and $T\geq0$ imply that $-Tb\leq Tx\leq Ta$.
Then, since $Y$ is $1$-max-normal, 
\[
\|Tx\|\leq\max\{\|Ta\|,\|Tb\|\}\leq\|T\|_{+}\max\{\|a\|,\|b\|\}\leq\|T\|_{+}(\|x\|+\varepsilon).
\]
Because $\varepsilon>0$ was chosen arbitrarily, we conclude that
$\|T\|_{+}\leq\|T\|\leq\|T\|_{+}$.

We prove the sufficiency of (2). Let $T\in B(X,Y)_{+}$. Let $x\in X$
and $\varepsilon>0$ be arbitrary, then there exists a $z\in X_{+}$
such that $\{-x,x\}\leq z$ and $\|z\|<\|x\|+\varepsilon$. Then,
since $T\geq0$, we see that $\{-Tx,Tx\}\leq Tz$, and therefore,
since $Y$ is absolutely monotone, 
\[
\|Tx\|\leq\|Tz\|\leq\|T\|_{+}\|z\|\leq\|T\|_{+}(\|x\|+\varepsilon).
\]
Because $\varepsilon>0$ was chosen arbitrarily, we conclude that
$\|T\|_{+}\leq\|T\|\leq\|T\|_{+}$.

We prove the sufficiency of (3). Let $x\in X$ be arbitrary. Since
$X$ is approximately $1$-sum-conormal, for every $\varepsilon>0$,
there exist $a,b\in X_{+}$ such that $x=a-b$ and $\|a\|+\|b\|<\|x\|+\varepsilon$.
For any $T\in B(X,Y)$, we have 
\[
\|Tx\|\leq\|Ta\|+\|Tb\|\leq\|T\|_{+}(\|a\|+\|b\|)\leq\|T\|_{+}(\|x\|+\varepsilon).
\]
Since $\varepsilon>0$ and $x\in X$ were chosen arbitrarily, we obtain
$\|T\|_{+}\leq\|T\|\leq\|T\|_{+}$. 
\end{proofof}

\section{Quasi-lattices and their basic properties\label{sec:Quasi-lattices-basic-properties}}

In this section we will define quasi-lattices, establish their basic
properties and provide a number of illustrative (non-)examples. 

Let $X$ be a pre-ordered Banach space and $A$ any subset of $X$.
For $x\in X$, by $A\leq x$ we mean that $a\leq x$ for all $a\in A$
and say \emph{$x$ is an upper bound of $A$. }We will use the Greek
letter `upsilon' to denote the set of all upper bounds of $A$, written
as  $\upsilon(A)$. If $x\in X$ is such that $A\leq x$ and, for
any $y\in X$, $A\leq y\leq x$ implies  $x=y$, we say that $x$\emph{
}is a \emph{minimal upper bound of $A$. }We will use the Greek letter
`mu' to denote the set of all minimal upper bounds of $A$, written
as $\mu(A)$. We note that $\upsilon(A)$ and $\mu(A)$ could be empty
for some $A\subseteq X$.

For any fixed $x,y\in X$, we define the function $\sigma_{x,y}:X\to\mathbb{R}_{\geq0}$
by $\sigma_{x,y}(z):=\|z-x\|+\|z-y\|$ for all $z\in X$, and note
that $\sigma_{x,y}(z)\geq\|x-y\|$ for all $x,y,z\in X$. We will
refer to $\sigma_{x,y}$ as \emph{the distance sum} \emph{to $x$
and $y$}. 

We introduce the following definitions and notation: 
\begin{defn}
\label{def:quasi-lattices}Let $X$ be a pre-ordered Banach space
with a closed cone.
\begin{enumerate}
\item We say that $X$ is an \emph{$\upsilon$-quasi-lattice }if, for every
pair of elements $x,y\in X$, $\upsilon(\{x,y\})$ is non-empty and
there exists a unique element $z\in\upsilon(\{x,y\})$ minimizing
$\sigma_{x,y}$ on $\upsilon(\{x,y\})$. The element $z$ will be
called the \emph{$\upsilon$-quasi-supremum }of $\{x,y\}$.
\item We say that $X$ is a\emph{ $\mu$-quasi-lattice }if, for every pair
of elements $x,y\in X$, $\mu(\{x,y\})$ is non-empty and there exists
a unique element $z\in\mu(\{x,y\})$ minimizing $\sigma_{x,y}$ on
$\mu(\{x,y\})$. The element $z$ will be called the \emph{$\mu$-quasi-supremum
}of $\{x,y\}$.
\end{enumerate}
\end{defn}
We immediately note that all Banach lattices are $\mu$-quasi-lattices:
\begin{prop}
\label{prop:Banach-lattices-are-mu-quasi}If $X$ is a lattice ordered
Banach space with a closed cone \textup{(}in particular, if $X$
is a Banach lattice\textup{)}, then $X$ is a $\mu$-quasi-lattice
and its lattice structure coincides with its $\mu$-quasi-lattice
structure.\end{prop}
\begin{proofof}
Since for every $x,y\in X$, $\mu(\{x,y\})=\{x\vee y\}$ is a singleton,
this is clear.\end{proofof}
\begin{rem}
If $X$ is a pre-ordered Banach space with a closed cone, then, for
$x,y\in X$, the set $\upsilon(\{x,y\})$ is closed and convex, and
hence techniques from convex optimization can be used to establish
whether a pre-ordered Banach space is an $\upsilon$-quasi-lattice
(cf. Theorem \ref{thm:reflexive-u-quasi-lattice}). The set $\mu(\{x,y\})$
need not be convex in general (cf.\,Example \ref{Example--comparability}),
and hence it is usually more difficult to determine whether or not
a space is a $\mu$-quasi-lattice than an $\upsilon$-quasi-lattice. 

Except in the case of monotone $\upsilon$-quasi-lattices which are
also $\mu$-quasi-lattices with coinciding $\upsilon$-- and $\mu$-quasi-lattice
structures (cf.\,Theorem \ref{thm:monotone-u-quasi-implies-mu-quasi}),
no further relationship is known between $\upsilon$-- and $\mu$-quasi-lattices.
Example \ref{example:mu-quasi-lattice-but-not-u-quasi-lattice} will
provide a Banach lattice, and hence $\mu$-quasi-lattice, that is
not an $\upsilon$-quasi-lattice. Furthermore, Example \ref{example:polynomial_cone}
will provide a non-monotone $\upsilon$-quasi-lattice, which exhibits
$\upsilon$-quasi-suprema that are not minimal, hence if this space
were a $\mu$-quasi-lattice (which is currently not known), then its
$\upsilon$-- and $\mu$-quasi-lattice structures will not coincide.
\end{rem}
To avoid repetition, we will often use the term quasi-lattice when
it is unimportant whether a space is an $\upsilon$-- or $\mu$-quasi-lattice,
i.e., a quasi-lattice is either an $\upsilon$-- or $\mu$-quasi-lattice.
In such cases we will refer to the relevant $\upsilon$-- or $\mu$-quasi-supremum
as just the quasi-supremum. When it is indeed important whether a
space is an\emph{ $\upsilon$}-- or $\mu$-quasi-lattice, we will
mention it explicitly.

The following notation will be used for both $\upsilon$-- and $\mu$-quasi-lattices.
Let $X$ be a quasi-lattice and $x,y\in X$ arbitrary. We will denote
the quasi-supremum of $\{x,y\}$ by $x\tilde{\vee}y$. This operation
is symmetric, i.e., $x\tilde{\vee}y=y\tilde{\vee}x$. We define the
\emph{quasi-infimum }of $\{x,y\}$ by $x\tilde{\wedge}y:=-((-x)\tilde{\vee}(-y))$.
It is elementary to see that $x\tilde{\wedge}y\le\{x,y\}$. We define
the \emph{quasi-absolute value} of $x$ by $\left\lceil x\right\rceil :=(-x)\tilde{\vee}(x)$.
We will often use the notation $x^{+}:=0\tilde{\vee}x$ and $x^{-}:=0\tilde{\vee}(-x)$. 

Before establishing the basic properties of quasi-lattices, we will
give a few examples of spaces that are (not) quasi-lattices.

The following is an example of a quasi-lattice that is not a Riesz
space, and hence not a Banach lattice:
\begin{exampleenv}
\label{Example-icecream}The space $\{\mathbb{R}^{3},\|\cdot\|_{2}\}$,
endowed with the Lorentz cone 
\[
C:=\{(x_{1},x_{2},x_{3}):x_{1}\geq(x_{2}^{2}+x_{3}^{2})^{1/2}\}.
\]
There are many minimal upper bounds of, e.g., $\{(0,0,0),(0,0,2)\}$
(cf.\,Example \ref{Example--comparability} and Proposition \ref{prop:Hilbert-sp-with-lorenz-not-banach-lattice}).
Hence no supremum exists, and this space is not a Riesz space. Another
method to establish this would be to note that $C$ has more than
distinct $3$ extreme rays, while every lattice cone in $\mathbb{R}^{3}$
has at most $3$ disinct extreme rays \cite[Theorem 1.45]{CnD}.

This space is (simultaneously an $\upsilon$-quasi-lattice and) a
$\mu$-quasi-lattice. Intuitively, this can be seen by taking arbitrary
elements, $x,y\in\mathbb{R}^{3}$, and seeing that there exists a
unique element in $\mu(\{x,y\})$ with least first coordinate, which
is then the quasi-supremum. It is possible give a more explicit proof,
but this is not needed in view of the general Theorem \ref{thm:hilbert-lorentz-collected-properties}
which is applicable to this example.
\end{exampleenv}
The following is an example of a Banach lattice, hence a $\mu$-quasi-lattice,
that is not an $\upsilon$-quasi-lattice:
\begin{exampleenv}
\label{example:mu-quasi-lattice-but-not-u-quasi-lattice}Consider
the space $\{\mathbb{R}^{3},\|\cdot\|_{\infty}\}$ with the standard
cone. Let $x:=(1,-1,0)$. Then 
\[
\upsilon(\{0,x\})=\{(z_{1},z_{2},z_{3})\in\mathbb{R}^{3}:z_{1}\geq1,\ z_{2},z_{3}\geq0\},
\]
and hence, for all $z\in\upsilon(\{0,x\})$, we see $\sigma_{0,x}(z)=\|z\|_{\infty}+\|(z_{1}-1,z_{2}+1,z_{3})\|_{\infty}\geq1+1=2$.
But, for every $t\in[0,1]$, $\{0,x\}\leq z_{t}:=(1,0,t)$ is such
that 
\[
\sigma_{x,0}(z_{t})=\|z_{t}-x\|_{\infty}+\|z_{t}-0\|_{\infty}=\|(0,1,t)\|_{\infty}+\|(1,0,t)\|_{\infty}=2,
\]
so that there exists no unique upper bound of $\{x,0\}$ minimizing
the distance sum to $x$ and $0$. We conclude that this space is
not an $\upsilon$-quasi-lattice.
\end{exampleenv}
There do exist ordered Banach spaces endowed with closed proper generating
cones that are not normed Riesz spaces, nor $\mu$-quasi-lattices
or $\upsilon$-quasi-lattices:
\begin{exampleenv}
\label{Non-example}Let $\{\mathbb{R}^{3},\|\cdot\|_{\infty}\}$ be
endowed with the cone defined by the four extreme rays $\{(\pm1,\pm1,1)\}$.
Let $x:=(0,0,0)$ and $y:=(2,0,0)$. It can be seen that $\mu(\{x,y\})=\{(1,t,1)\in\mathbb{R}^{3}:t\in[-1,1]\}$.
Since this set is not a singleton, this space is not a Riesz space.
Moreover, $\sigma_{x,y}$ takes the constant value $2$ on $\mu(\{x,y\})$,
and hence there does not exist a unique element minimizing $\sigma_{x,y}$
on $\mu(\{x,y\})$. Therefore this space is not a $\mu$-quasi-lattice.
Furthermore, if $z\in\upsilon(\{x,y\})$ and $z_{3}>1$, then $\sigma_{x,y}(z)>2$,
and since $\upsilon(\{x,y\})\cap\{z\in\mathbb{R}^{3}:z_{3}\leq1\}=\mu(\{x,y\})$,
all minimizers of $\sigma_{x,y}$ in $\upsilon(\{x,y\})$ must be
elements of $\mu(\{x,y\})$. Since $\sigma_{x,y}$ is constant on
$\mu(\{x,y\})$, this space is not an $\upsilon$-quasi-lattice.
\end{exampleenv}
The following results establish some basic properties of quasi-lattices.
\begin{prop}
\label{prop:quasi-lattice-cones-are-proper-and-generating}If $X$
is a quasi-lattice, then $X_{+}$ is a proper and generating cone.\end{prop}
\begin{proofof}
If $X_{+}$ is not proper, there exists an $x\in X$ such that $x>0$
and $-x>0$. Let $z\geq\{0,x\}$ be arbitrary. Then $z-x>z\geq0>x$,
so that $z>z-x\geq\{0,x\}$. Hence no upper bound of $\{0,x\}$ is
minimal and therefore $X$ cannot be a $\mu$-quasi-lattice.

Moreover, every $z\in\{\lambda x:\lambda\in[-1,1]\}$ minimizes $\sigma_{-x,x}$
on $\upsilon(\{x,-x\})$, therefore $X$ cannot be an $\upsilon$-quasi-lattice
either.

For all $x\in X$, since $\upsilon(\{x,0\})$ is non-empty, taking
any $z\in\upsilon(\{x,0\})$ and writing $x=z-(z-x)$ shows that $X_{+}$
is generating in $X$.
\end{proofof}
Surprisingly, many elementary Riesz space properties have direct analogues
in quasi-lattices. Many of the proofs below follow arguments from
\cite[Sections 5, 6]{ZaanenIntro} nearly verbatim.
\begin{thm}
\label{thm:formulas1}Let $X$ be a quasi-lattice, and $x,y,z\in X$
arbitrary. Then:
\begin{enumerate}
\item $x\tilde{\vee}x=x\tilde{\wedge}x=x$. 
\item For $\alpha\geq0$, $(\alpha x)\tilde{\vee}(\alpha y)=\alpha(x\tilde{\vee}y)$
and $(\alpha x)\tilde{\wedge}(\alpha y)=\alpha(x\tilde{\wedge}y)$.
\item For $\alpha\leq0$, $(\alpha x)\tilde{\vee}(\alpha y)=\alpha(x\tilde{\wedge}y)$
and $(\alpha x)\tilde{\wedge}(\alpha y)=\alpha(x\tilde{\vee}y)$.
\item $(x\tilde{\vee}y)+z=(x+z)\tilde{\vee}(y+z)$ and $(x\tilde{\wedge}y)+z=(x+z)\tilde{\wedge}(y+z)$.
\item $x^{\pm}\geq0$, $x^{-}=(-x)^{+}$.
\item $\left\lceil x\right\rceil \geq0$ and, for all $\alpha\in\mathbb{R}$,
$\left\lceil \alpha x\right\rceil =|\alpha|\left\lceil x\right\rceil $.
In particular $\left\lceil -x\right\rceil =\left\lceil x\right\rceil $.
\item $x=x^{+}-x^{-}$; $x^{+}\tilde{\wedge}x^{-}=0$ and $\left\lceil x\right\rceil =x^{+}+x^{-}$.
\item If $x\geq0$, then $x\tilde{\wedge}0=0$ and $x=x^{+}=\left\lceil x\right\rceil $.
\item $\left\lceil \left\lceil x\right\rceil \right\rceil =\left\lceil x\right\rceil $.
\item $x\tilde{\vee}y+x\tilde{\wedge}y=x+y$ and $x\tilde{\vee}y-x\tilde{\wedge}y=\left\lceil x-y\right\rceil $.
\item $x\tilde{\vee}y=\frac{1}{2}(x+y)+\frac{1}{2}\left\lceil x-y\right\rceil $
and $x\tilde{\wedge}y=\frac{1}{2}(x+y)-\frac{1}{2}\left\lceil x-y\right\rceil $.
\end{enumerate}
\end{thm}
\begin{proofof}
Assertion (1) follows from $x\leq x$ and the fact that $\sigma_{x,x}(x)=0$
and $\sigma_{x,x}(y)>0$ for all $y\neq x$.

We prove the assertion (2) for $\mu$-quasi-lattices. The case $\alpha=0$
follows from (1), hence we assume $\alpha>0$. By definition, $x\tilde{\vee}y$
is a minimal upper bound of $\{x,y\}$. Since $\alpha>0$, the element
$\alpha(x\tilde{\vee}y)$ is then a minimal upper bound of $\{\alpha x,\alpha y\}$.
Suppose that $\alpha(x\tilde{\vee}y)\neq(\alpha x)\tilde{\vee}(\alpha y)$,
then there exists a minimal upper bound of $\{\alpha x,\alpha y\}$,
say $z_{0}$, such that 
\[
\sigma_{\alpha x,\alpha y}(z_{0})=\|z_{0}-\alpha x\|+\|z_{0}-\alpha y\|<\|\alpha(x\tilde{\vee}y)-\alpha x\|+\|\alpha(x\tilde{\vee}y)-\alpha y\|.
\]
But then $\alpha^{-1}z_{0}$ is a minimal upper bound for $\{x,y\}$,
and 
\[
\sigma_{x,y}(\alpha^{-1}z_{0})=\|\alpha^{-1}z_{0}-x\|+\|\alpha^{-1}z_{0}-y\|<\|(x\tilde{\vee}y)-x\|+\|(x\tilde{\vee}y)-y\|,
\]
contradicting the definition of $x\tilde{\vee}y\in\mu(\{x,y\})$ as
the unique element minimizing $\sigma_{x,y}$ on $\mu(\{x,y\})$.
We conclude that $(\alpha x)\tilde{\vee}(\alpha y)=\alpha(x\tilde{\vee}y)$.
The same argument holds for $\upsilon$-quasi-lattices by ignoring
the word `minimal' in the previous argument. By using what was just
established, showing that $(\alpha x)\tilde{\wedge}(\alpha y)=\alpha(x\tilde{\wedge}y)$
holds is an elementary calculation.

The assertion (3) follows from applying (2) with $\beta:=-\alpha\geq0$.

The assertion (4) follows from the translation invariance of both
the metric defined by the norm and the partial order, and (5) is immediate
from the definitions. 

To establish (6), we notice that $\{x,-x\}\leq\left\lceil x\right\rceil $
implies $0\leq x-x\leq2\left\lceil x\right\rceil $. The second part
follows by noticing that $(-\alpha x)\tilde{\vee}(\alpha x)=(-|\alpha|x)\tilde{\vee}(|\alpha|x)$
and applying (2).

We prove (7). By (4), $x^{+}-x=(x\tilde{\vee}0)-x=(x-x)\tilde{\vee}(-x)=0\tilde{\vee}(-x)=x^{-}$,
so $x=x^{+}-x^{-}$. By this, we then have $0=-x^{-}+x^{-}=x\tilde{\wedge}0+x^{-}=(x+x^{-})\tilde{\wedge}(x^{-})=(x^{+})\tilde{\wedge}(x^{-})$.
By (2) and (4), $\left\lceil x\right\rceil =(-x)\tilde{\vee}x=(-x)\tilde{\vee}x+x-x=0\tilde{\vee}(2x)-x=2x^{+}-x^{+}+x^{-}=x^{+}+x^{-}.$

We prove (8). Let $x\geq0$, then $0$ is an upper bound of $\{0,-x\}$.
Moreover, since the cone is proper, $0$ is a minimal upper bound
for $\{0,-x\}$. But, for any $z\in X$ (and in particular all (minimal)
upper bounds of $\{0,-x\}$), we have 
\[
\sigma_{-x,0}(0)=\|0-(-x)\|+\|0-0\|=\|0-(-x)\|\leq\sigma_{-x,0}(z).
\]
Hence we have $0=(-x)\tilde{\vee}0=x^{-}$, and hence, by (7), $x=x^{+}=\left\lceil x\right\rceil $.

The assertion (9) follows from (6) and (8).

We prove (10). We observe that $x\tilde{\vee}y=((x-y)\tilde{\vee}0)+y=(x-y)^{+}+y$,
and $x\tilde{\wedge}y=x+(0\tilde{\wedge}(y-x))=x-(x-y)^{+}$. Adding
these two equations yields $x\tilde{\vee}y+x\tilde{\wedge}y=x+y$,
and subtracting gives $x\tilde{\vee}y-x\tilde{\wedge}y=2(x-y)^{+}+y-x=(2(x-y)\tilde{\vee}0)-(x-y)=((x-y)\tilde{\vee}(-(x-y))=\left\lceil x-y\right\rceil $. 

The assertion (11) follows by adding and subtracting the equations
established in (10).
\end{proofof}
In a sense the more interesting results concerning quasi-lattices
are ones outlining how they may differ from Riesz spaces and Banach
lattices. An important remark, that may at first sight be counterintuitive,
is the following: For elements $x,y,z$ in a quasi-lattice, $x\leq z$
and $y\leq z$ does not, in general, imply that $x\tilde{\vee}y\leq z$.
The following example shows how this may happen:
\begin{exampleenv}
\label{Example--comparability}We continue with Example \ref{Example-icecream}.
Let $x=(0,0,0)$ and $y=(0,0,2)$, then $x\tilde{\vee}y=(1,0,1)$.
The set of minimal upper bounds of $\{x,y\}$ forms a branch of a
hyperbola. Choosing $z$ from this hyperbola such that $z$ and $x\tilde{\vee}y$
are not comparable, say any $z=(\sqrt{t^{2}+1},\pm t,1)$ with $t>0$,
we see that, although $x\leq z$ and $y\leq z$, it does not hold
that $x\tilde{\vee}y\leq z$.
\end{exampleenv}
The previous example shows how it may sometimes happen in quasi-lattices
that the quasi-supremum operation is not monotone: $x\leq y$ does
not necessarily imply $x^{+}\leq y^{+}$. We can therefore not expect
distributive laws, Birkhoff type inequalities or the Riesz decomposition
property to hold in general quasi-lattices.

The following example shows how a quasi-supremum operation need not
even be associative:
\begin{exampleenv}
Let $\{\mathbb{R}^{3},\|\cdot\|_{2}\}$ be endowed with a `half Lorentz
cone' 
\[
C:=\{(x_{1},x_{2},x_{3}):x_{1}\geq(x_{2}^{2}+x_{3}^{2})^{1/2},x_{2}\geq0\}.
\]
 By Corollary \ref{cor:monotone-reflexive-is-a-mu-quasi-lattice},
this space is a $\mu$-quasi-lattice. 

For any $x,y\in\mathbb{R}^{3}$, we claim that $(x\tilde{\vee}y)_{2}=\max\{x_{2},y_{2}\}$.
To this end, let $z\geq\{x,y\}$ be arbitrary and define $z':=(z_{1},\max\{x_{2},y_{2}\},z_{3})$.
We first show that $z'\geq\{x,y\}$. Firstly, $z_{2}'-x_{2}=\max\{x_{2},y_{2}\}-x_{2}\geq0$
and $z_{2}'-y_{2}=\max\{x_{2},y_{2}\}-y_{2}\geq0$. Since $z_{2}-x_{2}\geq0$
and $z_{2}-y_{2}\geq0$, we also have $z_{2}\geq z_{2}'$. Also, where
we use the fact that $(z_{2}-z_{2}')(z_{2}-x_{2})\geq0$ and $(z_{2}-z_{2}')^2\geq0$
in the last step, 
\begin{eqnarray*}
z_{1}'-x_{1} & = & z_{1}-x_{1}\\
 & \geq & \sqrt{(z_{2}-x_{2})^{2}+(z_{3}-x_{3})^{2}}\\
 & = & \sqrt{(z_{2}-z_{2}'+z_{2}'-x_{2})^{2}+(z_{3}-x_{3})^{2}}\\
 & = & \sqrt{(z_{2}-z_{2}')^{2}+2(z_{2}-z_{2}')(z_{2}-x_{2})+(z_{2}'-x_{2})^{2}+(z_{3}'-x_{3})^{2}}\\
 & \geq & \sqrt{(z_{2}'-x_{2})^{2}+(z_{3}'-x_{3})^{2}}.
\end{eqnarray*}
Similarly $z_{1}'-y_{1}\geq\sqrt{(z_{2}'-y_{2})^{2}+(z_{3}'-y_{3})^{2}}$,
so that $z'\geq\{x,y\}$. We claim that $\sigma_{x,y}(z)\geq\sigma_{x,y}(z')$.
Indeed, again since $(z_{2}-z_{2}')(z_{2}-x_{2})\geq0$ and $(z_{2}-z_{2}')^2\geq0$,
\begin{eqnarray*}
 &  & \|z-x\|_{2}\\
 & = & \sqrt{(z_{1}-x_{1})^{2}+(z_{2}-x_{2})^{2}+(z_{3}-x_{2})^{2}}\\
 & = & \sqrt{(z_{1}-x_{1})^{2}+(z_{2}-z_{2}'+z_{2}'-x_{2})^{2}+(z_{3}-x_{2})^{2}}\\
 & = & \sqrt{(z_{1}-x_{1})^{2}+(z_{2}-z_{2}')^{2}+2(z_{2}-z_{2}')(z_{2}'-x_{2})+(z_{2}'-x_{2})^{2}+(z_{3}-x_{2})^{2}}\\
 & \geq & \sqrt{(z_{1}-x_{1})^{2}+(z_{2}'-x_{2})^{2}+(z_{3}-x_{2})^{2}}\\
 & = & \sqrt{(z_{1}'-x_{1})^{2}+(z_{2}'-x_{2})^{2}+(z_{3}'-x_{2})^{2}}\\
 & = & \|z'-x\|_{2}.
\end{eqnarray*}
Similarly we have $\|z-y\|_{2}\geq\|z'-y\|_{2}$. Therefore $\sigma_{x,y}(z)=\|z-x\|_{2}+\|z-y\|_{2}\geq\sigma_{x,y}(z')$.
We conclude that $(x\tilde{\vee}y)_{2}=\max\{x_{2},y_{2}\}$, else, by the above construction, 
there would exist an upper bound of $\{x,y\}$ different from $x\tilde{\vee}y$, but which also minimizes $\sigma_{x,y}$ on $\mu(\{x,y\})$. 

Now let $a:=(0,0,0)$, $b:=(0,-1,1)$ and $c:=(0-1,-1)$. Using what
was just proven and the fact that the space is a $\mu$-quasi-lattice,
it can be seen that $a\tilde{\vee}b$ must be an element of the plane $\{x\in\mathbb{R}^{3}:x_{2}=0\}$
and must be a minimal upper bound of $\{a,b\}$. The minimal upper
bounds of $\{a,b\}$ that are elements of $\{x\in\mathbb{R}^{3}:x_{2}=0\}$
can be parameterized by $\gamma:t\mapsto(\sqrt{1+(1-t)^{2}},0,t)$
with $t\in(-\infty,1]$ and the function $t\mapsto\sigma_{a,b}(\gamma(t))$
attains its minimum at $t=\sqrt{3}-1$. Therefore $a\tilde{\vee}b=(2\sqrt{2-\sqrt{3}},0,\sqrt{3}-1)$.
Again, using similar reasoning, it can be verified (using a computer
algebra system) that $(a\tilde{\vee}b)\tilde{\vee}c=(\sqrt{1+\left(1+\kappa\right)^{2}},0,\kappa)$,
where $\kappa:=23^{-1}\left(-29-8\sqrt{2}+9\sqrt{3}+12\sqrt{6}\right)$.
Also, since $(1,-1,0)$ is the only minimal upper bound of $\{b,c\}$
that is an element of the plane $\{x\in\mathbb{R}^{3}:x_{2}=-1\}$,
we must have $b\tilde{\vee}c=(1,-1,0)$. It can then be verified that
$a\tilde{\vee}(b\tilde{\vee}c)=(2,0,0)$. We conclude that $a\tilde{\vee}(b\tilde{\vee}c)\neq(a\tilde{\vee}b)\tilde{\vee}c$.

\end{exampleenv}
The triangle and reverse triangle inequality take the following form
in quasi-lattices. They reduce to the familiar ones in lattice-ordered
$\mu$-quasi-lattices.
\begin{thm}
\textup{(}Triangle and reverse triangle inequality\textup{)} Let
$X$ be a quasi-lattice and $x,y\in X$ be arbitrary. Then 
\[
\{x+y,-(x+y)\}\leq\left\lceil x\right\rceil +\left\lceil y\right\rceil ,
\]
and
\[
\{x-\left\lceil y\right\rceil ,-x-\left\lceil y\right\rceil ,y-\left\lceil x\right\rceil ,-y-\left\lceil x\right\rceil \}\leq\left\lceil x\pm y\right\rceil .
\]
\end{thm}
\begin{proofof}
By Theorem \ref{thm:formulas1} (7), for all $z\in X$, we have $\left\lceil z\right\rceil \geq z^{\pm}\geq\pm z$,
and hence we obtain $\left\lceil x\right\rceil +\left\lceil y\right\rceil \geq x^{+}+y^{+}\geq x+y$
and $\left\lceil x\right\rceil +\left\lceil y\right\rceil \geq x^{-}+y^{-}\geq-x-y$.
Therefore $\left\lceil x\right\rceil +\left\lceil y\right\rceil $
is an upper bound of $\{x+y,-(x+y)\}$.

To establish the second inequality, we use what was just established
to see, by Theorem \ref{thm:formulas1} (6) and (9), that $\{x,-x\}=\{(x\pm y)\mp y,-((x\pm y)\mp y)\}\leq\left\lceil x\pm y\right\rceil +\left\lceil \mp y\right\rceil =\left\lceil x\pm y\right\rceil +\left\lceil y\right\rceil $.
Hence $\{x-\left\lceil y\right\rceil ,-x-\left\lceil y\right\rceil \}\leq\left\lceil x\pm y\right\rceil $.
Similarly, $\{y-\left\lceil x\right\rceil ,-y-\left\lceil x\right\rceil \}\leq\left\lceil x\pm y\right\rceil $,
and finally we conclude that $\{x-\left\lceil y\right\rceil ,-x-\left\lceil y\right\rceil ,y-\left\lceil x\right\rceil ,-y-\left\lceil x\right\rceil \}\leq\left\lceil x\pm y\right\rceil $.
\end{proofof}
The following result allows us to conclude that monotone $\upsilon$-quasi-lattices
are in fact $\mu$-quasi-lattices:
\begin{thm}
\label{thm:monotone-u-quasi-implies-mu-quasi}Every monotone $\upsilon$-quasi-lattice
is a $\mu$-quasi-lattice, and its $\upsilon$-- and $\mu$-quasi-lattice
structures coincide.\end{thm}
\begin{proofof}
We first claim that, if $X$ is a monotone ordered Banach space, then,
for $x,y\in X$, if $z_{0}\in\upsilon(\{x,y\})$ is such that $\|z-x\|+\|z-y\|>\|z_{0}-x\|+\|z_{0}-y\|$
for all $z\in\upsilon(\{x,y\})$ with $z\neq z_{0}$, then $z_{0}$
is a minimal upper bound of $\{x,y\}$.

As to this, by translating, we may assume that $y=0$. Let $z\in X$
be any element satisfying $\{x,0\}\leq z\leq z_{0}$. Then $0\leq z\leq z_{0}$
and $0\leq z-x\leq z_{0}-x$. By monotonicity, $\|z\|\leq\|z_{0}\|$
and $\|z-x\|\leq\|z_{0}-x\|$, so that $\|z\|+\|z-x\|\leq\|z_{0}\|+\|z_{0}-x\|$.
The hypothesis on $z_{0}$ then implies that $z=z_{0}$. Hence $z_{0}$
is a minimal upper bound of $\{x,y\}$, establishing the claim.

Let $X$ be a monotone $\upsilon$-quasi-lattice and $x,y\in X$ arbitrary.
By the above claim, the $\upsilon$-quasi-supremum of $\{x,y\}$ is
a minimal upper bound of $\{x,y\}$. Since $\mu(\{x,y\})\subseteq\upsilon(\{x,y\})$
we have that the $\upsilon$-quasi-supremum of $\{x,y\}$ is also
the $\mu$-quasi-supremum. We conclude that $X$ is also a $\mu$-quasi-lattice,
and that its $\upsilon$-- and $\mu$-quasi-lattice structures coincide.
\end{proofof}
The following example shows that there exist $\upsilon$-quasi-lattices
in which some $\upsilon$-quasi-suprema are not minimal upper bounds.
\begin{exampleenv}
\label{example:polynomial_cone}Consider the space $\{\mathbb{R}^{3},\|\cdot\|_{2}\}$,
endowed with the cone 
\[
C:=\{(a,b,c)\in\mathbb{R}^{3}:ax^{2}+bx+c\geq0\mbox{ for all }x\in[0,1]\}.
\]
By Theorem \ref{thm:reflexive-u-quasi-lattice}, this space is an
$\upsilon$-quasi-lattice. Let $x:=(0,1,0)$, $y:=(0,-1,1)$. It can
be verified (using a computer algebra system) that 
\[
x\tilde{\vee}y=(2^{-1}(2-\sqrt{3}),-2^{-1}(2-\sqrt{3}),1),
\]
 while $\{x,y\}\leq(1,-1,1)<x\tilde{\vee}y$. Therefore $x\tilde{\vee}y\notin\mu(\{x,y\})$.

By comparing the norms of the elements in $0\leq(0,-1,1)\leq(0,0,1)$,
we see that this space is not monotone. We can therefore not draw
any  conclusion from Theorem \ref{thm:monotone-u-quasi-implies-mu-quasi}
as to whether this space is a $\mu$-quasi-lattice. A valid conclusion
we may draw is that, if this space is indeed also a $\mu$-quasi-lattice
in addition to being an $\upsilon$-quasi-lattice, its $\mu$-- and
$\upsilon$-quasi-lattice structures will \emph{not} coincide.
\end{exampleenv}

\section{A concrete class of quasi-lattices\label{sec:More-on-quasi-lattices}}

In the previous section we have already noted that lattice ordered
Banach spaces are $\mu$-quasi-lattices (cf.\,Proposition \ref{prop:Banach-lattices-are-mu-quasi})
and gave a number of examples of quasi-lattices. We begin this section
by proving that quite a large class of (not necessarily lattice ordered)
ordered Banach spaces with closed generating cones are in fact quasi-lattices.
Afterwards, we briefly investigate conditions under which a space
has a quasi-lattice as a dual, or is the dual of a quasi-lattice.

We recall that a normed space $X$ is \emph{strictly convex }or \emph{rotund}
if, for $x,y\in X$, $\|x+y\|=\|x\|+\|y\|$ implies that either $x$
or $y$ is a non-negative multiple of the other \cite[Definition 5.1.1, Proposition 5.1.11]{Megginson}. 

The following theorem shows that there exist relatively many quasi-lattices:
\begin{thm}
\label{thm:reflexive-u-quasi-lattice}Every strictly convex reflexive
ordered Banach space $X$ with a closed proper generating cone is
an $\upsilon$-quasi-lattice.\end{thm}
\begin{proofof}
We need to prove that every pair of elements $x_{0},y_{0}\in X$ has
an $\upsilon$-quasi-supremum in $X$. 

If $x_{0}$ and $y_{0}$ are comparable,  by exchanging the roles
of $x_{0}$ and $y_{0}$ if necessary, we may assume $x_{0}\leq y_{0}$.
We may further assume that $x_{0}=0$ by translating over $-x_{0}$.
We will denote the distance sum to $0$ and $y_{0}$ by $\sigma$
instead of $\sigma_{0,y_{0}}$.

If $y_{0}=0$, then $\sigma(z)=0$ if and only if $z=0$, so that
$0\tilde{\vee}0=0$. If $0\neq y_{0}\geq0$, we have that $y_{0}\in\upsilon(\{0,y_{0}\})$
and, for all $z\in\upsilon(\{0,y_{0}\})$, we have $\sigma(z)=\|y_{0}-z\|+\|z\|\geq\|y_{0}\|=\sigma(y_{0})$.
Suppose that $z_{0}\in\upsilon(\{0,y_{0}\})$ is such that $\sigma(y_{0})=\sigma(z_{0})$.
We must have $z_{0}\neq0$, else $0\leq y_{0}\leq z_{0}=0$ hence,
since $X_{+}$ is proper, $y_{0}=0$, while $y_{0}\neq0$. Then, since
\[
\|y_{0}-z_{0}+z_{0}\|=\|y_{0}\|=\sigma(y_{0})=\sigma(z_{0})=\|y_{0}-z_{0}\|+\|z_{0}\|,
\]
by strict convexity we obtain $y_{0}-z_{0}=\lambda z_{0}$ for some
$\lambda\geq0$. Hence $z_{0}\geq y_{0}=(1+\lambda)z_{0}\geq z_{0}$
and then, since $X_{+}$ is proper, $y_{0}=z_{0}$. Therefore $0\tilde{\vee}y_{0}=y_{0}$.

We consider the case where neither $x_{0}\leq y_{0}$ nor $y_{0}\leq x_{0}$.
Again, by translating, we may assume without loss of generality that
$x_{0}=0$, and that neither $y_{0}\leq0$ nor $0\leq y_{0}$. We
again denote the distance sum to $0$ and $y_{0}$ by $\sigma$ instead
of $\sigma_{0,y_{0}}$.

Since $X_{+}$ is generating, $\upsilon(\{y_{0},0\})=X_{+}\cap(y_{0}+X_{+})$
is non-empty, hence let $z_{0}\in X_{+}\cap(y+X_{+})$. Consider the
non-empty closed bounded and convex set 
\[
K:=X_{+}\cap(y_{0}+X_{+})\cap\{x\in X:\sigma(x)\leq\sigma(z_{0})\}.
\]
We note that $0,y_{0}\notin K$, since we had assumed that neither
$y_{0}\leq0$ nor $0\leq y_{0}$ holds. 

The function $\sigma$ is continuous and convex and, since $K$ is
bounded closed and convex and $X$ is reflexive, by \cite[Theorem 2.11]{convexity-n-optimization},
there exists an element $z_{m}\in K$ minimizing $\sigma$ on $K$.
We claim that $z_{m}$ is the unique minimizer of $\sigma$ on $K$.
To prove this claim it is sufficient to establish that $\sigma$ is
strictly convex on $K$, i.e., if $z,z'\in K$ with $z\neq z'$ and
$t\in(0,1)$, then $\sigma(tz+(1-t)z')<t\sigma(z)+(1-t)(z')$.

We first claim that the line $\mathbb{R}y_{0}$ does not intersect
$K$. Indeed, if $\lambda y_{0}\in K$ for some $\lambda\in\mathbb{R}$,
then we must have $\lambda\neq0$, since $0\notin K$. But then $\lambda y_{0}\in K\subseteq X_{+}$
implies that either $y_{0}\leq0$ or $0\leq y_{0}$, contrary to our
assumption that neither $y_{0}\leq0$ nor $0\leq y_{0}$. 

We now prove that $\sigma$ is strictly convex on $K$. Let $z,z'\in K$
be arbitrary but distinct and $t\in(0,1)$. If $z\neq\lambda z'$
for all $\lambda\geq0$, then, by strict convexity of $X$, $\|tz+(1-t)z'\|<t\|z\|+(1-t)\|z'\|$,
and hence $\sigma(tz+(1-t)z')<t\sigma(z)+(1-t)\sigma(z')$. On the
other hand, if $z'=\lambda z$ for some $\lambda\geq0$, we must have
that $\lambda\neq1$ (since $z\neq z'$) and $\lambda\neq0$ (since
$0\notin K$). Therefore, supposing that 
\[
\|(1-t)(y_{0}-z)+t(y_{0}-z')\|=(1-t)\|y_{0}-z\|+t\|y_{0}-z'\|,
\]
 by strict convexity of $X$, we obtain $(1-t)(y_{0}-z)=\rho t(y_{0}-z')$
for some $\rho>0$ (if $\rho=0$, then $y_{0}=z\in K$ contradicts
$y_{0}\notin K$). By rewriting, we obtain $((1-t)-\rho t)y_{0}=((1-t)-\rho t\lambda)z$.
If $((1-t)-\rho t\lambda)=0$, then $((1-t)-\rho t)\neq0$ since $\lambda\neq1$
and $\rho t\neq0$, and hence $y_{0}=0$, contradicting the assumption
that neither $y_{0}\leq0$ nor $0\leq y_{0}$. Therefore $((1-t)-\rho t\lambda)\neq0$,
and $z\in K\cap\mathbb{R}y_{0}$, contracting the fact that $K$ and
$\mathbb{R}y_{0}$ are disjoint. Therefore, we must have $\|(1-t)y_{0}-(1-t)z+ty_{0}-tz'\|<(1-t)\|y_{0}-z\|+t\|y_{0}-z'\|$,
and hence $\sigma(tz+(1-t)z')<t\sigma(z)+(1-t)\sigma(z')$. 

We conclude that $\sigma$ is strictly convex on $K$, and that $z_{m}\in K$
is the unique minimizer of $\sigma$ on $K$. Then clearly $z_{m}$
is also the unique minimizer of $\sigma$ on $\upsilon(\{0,y_{0}\})$.
\end{proofof}
Theorem \ref{thm:reflexive-u-quasi-lattice} and Theorem \ref{thm:monotone-u-quasi-implies-mu-quasi}
together yield the following two corollaries:
\begin{cor}
\label{cor:monotone-reflexive-is-a-mu-quasi-lattice}Every strictly
convex reflexive monotone ordered Banach space with a closed proper
generating cone is both an $\upsilon$-quasi-lattice and a $\mu$-quasi-lattice
\textup{(}and its $\upsilon$-- and $\mu$-quasi-lattice structures
coincide\textup{)}.
\end{cor}
{}
\begin{cor}
For $1<p<\infty$, every $L^{p}$-space endowed with a closed proper
generating cone is an $\upsilon$-quasi-lattice. In particular, every
$\ell^{p}$-space and every space $\{\mathbb{R}^{n},\|\cdot\|_{p}\}$
that is endowed with a closed proper generating cone is an $\upsilon$-quasi-lattice.
If, in addition, the space is monotone, it is also a $\mu$-quasi-lattice
\textup{(}and its $\upsilon$-- and $\mu$-quasi-lattice structures
coincide\textup{)}.\end{cor}
\begin{proofof}
That an $L^{p}$-space is strictly convex for every $1<p<\infty$
is a consequence of \cite[Theorem 5.2.11]{Megginson}. The result
then follows from the previous theorem and Theorem \ref{thm:monotone-u-quasi-implies-mu-quasi}.
\end{proofof}
The remainder of this section will be devoted to dual considerations,
specifically to the question of when the dual of a pre-ordered Banach
space is a quasi-lattice. The following result gives necessary conditions
for this to be the case.
\begin{prop}
\label{prop:dual-pseudolattice-implies-normality}If a pre-ordered
Banach space $X$ with a closed cone has a quasi-lattice as dual,
then: 
\begin{enumerate}
\item There exists an $\alpha>0$ such that $X$ is $\alpha$-max-normal.
\item $X_{+}-X_{+}$ is dense in $X$.
\end{enumerate}
\end{prop}
\begin{proofof}
By Proposition \ref{prop:quasi-lattice-cones-are-proper-and-generating},
the dual cone is proper and generating. By part (1)(e) of Theorem
\ref{thm:normality-conormality-duality-properties}, there exists
an $\alpha>0$ such that $X$ is $\alpha$-max-normal. By \cite[Theorem 2.13(2)]{CnD},
$X_{+}-X_{+}$ is weakly dense in $X$. Since $X_{+}-X_{+}$ is convex,
its weak closure and norm closure coincide, and $X_{+}-X_{+}$ is
therefore norm dense in $X$.\end{proofof}
\begin{cor}
Let $X$ be a pre-ordered Banach space with a closed generating cone.
If $X$ has a quasi-lattice as dual, then there exists an $\alpha>0$
such that $X$ is $\alpha$-Ellis-Grosberg-Krein regular. \end{cor}
\begin{proofof}
By the previous result, there exists a $\beta>0$ such that $X$ is
$\beta$-max-normal. The cone $X_{+}$ was assumed to be generating,
and therefore $X$ is And\^o-regular. By Proposition \ref{prop:regularity-relationships}
there exists an $\alpha>0$ such that $X$ is $\alpha$-Ellis-Grosberg-Krein
regular.
\end{proofof}
The following theorem provides sufficient conditions for a pre-ordered
Banach space to have a quasi-lattice as dual. 

We recall that a normed space $X$ is \emph{smooth} if, for every
$x\in X$ with $\|x\|=1$, there exists a unique element $\phi\in X'$
with $\|\phi\|=1$ such that $\phi(x)=1$ \cite[Definition 5.4.1, Corollary 5.4.3]{Megginson}.
\begin{thm}
\label{thm:spaces__with_pseudo(ersatz)-lattices_as_dual}If, for some
$\alpha>0$, $X$ is an $\alpha$-normal smooth reflexive pre-ordered
Banach space with a closed cone such that $X_{+}-X_{+}$ is dense
in $X$, then its dual is an $\upsilon$-quasi-lattice.

If, in addition, $X$ is approximately $1$-conormal, its dual is
a $\mu$-quasi-lattice \textup{(}and its $\upsilon$-- and $\mu$-quasi-lattice
structures coincide\textup{)}.\end{thm}
\begin{proofof}
By \cite[Corollary 2.14, Theorem 2.40]{CnD}, the dual cone is proper
and generating in $X'$. That the dual cone is closed is elementary.
By \cite[Proposition 5.4.7]{Megginson}, $X'$ is strictly convex,
since $X$ was assumed to be smooth. Therefore $X'$ satisfies the
hypotheses of Theorem \ref{thm:reflexive-u-quasi-lattice}, and is
an $\upsilon$-quasi-lattice.

If we make the extra assumption that $X$ is approximately $1$-conormal,
then by part (2)(d) of Theorem \ref{thm:normality-conormality-duality-properties},
$X'$ is monotone. Then, by Theorem \ref{thm:monotone-u-quasi-implies-mu-quasi},
$X'$ is a $\mu$-quasi-lattice.
\end{proofof}
The following example shows that there exist $\upsilon$-quasi-lattices
that are not $\alpha$-normal for any $\alpha>0$. It cannot have
a quasi-lattice as dual, nor is it the dual of a quasi-lattice. Indeed,
by Proposition \ref{prop:dual-pseudolattice-implies-normality} its
dual is not a quasi-lattice. Moreover, by part (2)(e) of Theorem \ref{thm:normality-conormality-duality-properties},
since the space is not $\alpha$-normal for any $\alpha>0$, its cone
is not the dual cone of a pre-ordered Banach space with a closed generating
cone, and in particular, it is not the dual of a quasi-lattice.
\begin{exampleenv}
\label{exa:non-normal}Consider the following subset of $\ell^{2}$:
\[
C:=\left\{ x\in\ell^{2}:x_{1}\geq\left(\sum_{m=2}^{\infty}\frac{1}{m}x_{m}^{2}\right)^{1/2}\right\} .
\]

Clearly, $C\cap(-C)=\{0\}$ and $\lambda C\subseteq C$ for all $\lambda\geq0$.
Also, by Minkowski's inequality, $C+C\subseteq C$ so that we may
conclude that $C$ is a proper cone. For any $x\in\ell^{2}$, taking
$y:=\lambda(1,0,0,\ldots)\in C$ with $\lambda\geq|x_{1}|+\left(\sum_{m=2}^{\infty}\frac{1}{m}x_{m}^{2}\right)^{1/2}$,
we see $x=y-(y-x)\in C-C$, so that $C$ is generating in $\ell^{2}$.
Since the map $\rho_{0}:\ell^{2}\to\mathbb{R}$ defined by $\rho_{0}:x\mapsto\left(\sum_{m=2}^{\infty}\frac{1}{m}x_{m}^{2}\right)^{1/2}$
is a continuous seminorm, the map $\rho:x\mapsto x_{1}-\rho_{0}(x)$
is also continuous. Since $C=\rho^{-1}(\mathbb{R}_{\geq0})$, we conclude
that $C$ is closed. By Theorem \ref{thm:reflexive-u-quasi-lattice},
this space is an $\upsilon$-quasi-lattice.

We claim that this space is not $\alpha$-normal for any $\alpha>0$.
It is sufficient to show, for every $\alpha\geq1$, that there exist
$x,y\in\ell^{2}$ with $0\leq x\leq y$, such that $\|x\|>\alpha\|y\|$.
To this end, we set $y:=(2,0,\ldots)$. We define $x$ as follows:
let $\mathbb{N}\ni n_{\alpha}>(2\alpha)^{2}$ and $x=(1,0,\ldots,0,\sqrt{n_{\alpha}},0,\ldots)$
with $\sqrt{n_{\alpha}}$ occurring at the $n_{\alpha}$-th coordinate.
We then see that $0\leq x\leq y$, while 
\[
\|x\|=\left(1+n_{\alpha}\right)^{\frac{1}{2}}>n_{\alpha}^{\frac{1}{2}}>2\alpha=\alpha\|y\|.
\]
We conclude that this space is not $\alpha$-normal for any $\alpha>0$.
\end{exampleenv}

\section{A class of
quasi-lattices with absolutely monotone spaces of operators
\label{sec:Quasi-lattices-with-positively-attained-operator-norms}}

In this section we show that a real Hilbert space $\mathcal{H}$ endowed
with a Lorentz cone (defined below) is a $1$-absolutely Davies-Ng
regular $\mu$-quasi-lattice (that is not a Banach lattice if $\dim\mathcal{H}\geq3$).
Through an application of Theorem \ref{thm-(co)normality-conditions-onXandY-implies-normality-ir-B(X,Y)},
this will resolve the question posed in the introduction of whether
there exist non-Banach lattice pre-ordered Banach spaces $X$ and
$Y$ for which $B(X,Y)$ is absolutely monotone. 

Results established in this section will be collected in Theorem \ref{thm:hilbert-lorentz-collected-properties}.
In particular it will be shown that $\|x\|=\|\left\lceil x\right\rceil \|$
for all $x\in\mathcal{H}$ (which is analogous to the identity $\|x\|=\||x|\|$
which holds for all elements $x$ of a Banach lattice). Then, for
$\alpha>0$ and pre-ordered Banach spaces $X$ and $Y$ that are respectively
approximately $\alpha$-absolutely conormal and $\alpha$-absolutely
normal, the spaces of operators $B(X,\mathcal{H})$ and $B(\mathcal{H},Y)$
are shown to be $\alpha$-absolutely normal. Furthermore, if $\alpha=1$
(in particular if $X$ and $Y$ are Hilbert spaces endowed with Lorentz
cones), then the operator norms of $B(X,\mathcal{H})$ and $B(\mathcal{H},Y)$
are positively attained. 

We begin with the following lemma which outlines sufficient conditions
for establishing $1$-absolute Davies-Ng regularity of absolutely
monotone quasi-lattices:
\begin{lem}
\label{lem:sufficient-conditions-for-quasilattice-positively-attained-operator-norm}Let
$X$ be an absolutely monotone quasi-lattice satisfying $\|x\|=\|\left\lceil x\right\rceil \|$
for all $x\in X$. Then $X$ is $1$-absolutely Davies-Ng regular.\end{lem}
\begin{proofof}
The fact that $\|x\|=\|\left\lceil x\right\rceil \|$ for all $x\in X$
implies that $X$ is $1$-absolutely conormal. Therefore $X$ is $1$-absolutely
Davies-Ng regular. 
\end{proofof}
Every Banach lattice satisfies the hypotheses of the previous proposition.
The rest of this section will be devoted to proving that there exist
quasi-lattices that are not Banach lattices, but still satisfy the
hypothesis of the previous proposition.
\begin{defn}
Let $\mathcal{H}$ be a real Hilbert space. For a norm-one element
$v\in\mathcal{H}$, let $P$ be the orthogonal projection onto $\{v\}^{\bot}$.
We define the \emph{Lorentz cone }
\[
\mathcal{L}_{v}:=\{x\in\mathcal{H}:\langle v|x\rangle\geq\|Px\|\}.
\]
As in Example \ref{exa:non-normal}, it is elementary to see that
this cone is closed, proper and generating in $\mathcal{H}$. 
\end{defn}
It is widely known that the Hilbert space $\mathbb{R}^{3}$ ordered
by the Lorentz cone $\mathcal{L}_{e_{1}}\subseteq\mathbb{R}^{3}$
is not a Riesz space (cf.\,Example \ref{Example-icecream}). This
is actually true for arbitrary Hilbert spaces endowed with a Lorentz
cone as we will now proceed to show. The following two lemmas will
be used in the proof of Proposition \ref{prop:Hilbert-sp-with-lorenz-not-banach-lattice}
which establishes this fact.
\begin{lem}
Let $\mathcal{H}$ be a real Hilbert space endowed with a Lorentz
cone $\mathcal{L}_{v}$ where $v\in\mathcal{H}$ is such that $\|v\|=1$.
If $x\in\mathcal{L}_{v}$ is such that $\langle x|v\rangle=\|Px\|$
and $z_{1},z_{2}\in\mathcal{L}_{v}$ are such that $x=z_{1}+z_{2}$,
then $z_{1},z_{2}\in\{\lambda x:\lambda\geq0\}$.\end{lem}
\begin{proofof}
Let $P$ be the orthogonal projection onto $\{v\}^{\bot}$. If $x=0$,
since $\mathcal{L}_{v}$ is proper, the statement is clear. Let $0\neq x\in\mathcal{L}_{v}$
be such that $\langle x|v\rangle=\|Px\|$. Then $\langle x|v\rangle=\|Px\|>0$,
else $x=0$. Suppose $z_{1,}z_{2}\in\mathcal{L}_{v}$ are such that
$x=z_{1}+z_{2}$. Then
\[
\langle x|v\rangle=\langle z_{1}+z_{2}|v\rangle\geq\|Pz_{1}\|+\|Pz_{2}\|\geq\|P(z_{1}+z_{2})\|=\|Px\|=\langle x|v\rangle.
\]
Therefore $\|Pz_{1}\|+\|Pz_{2}\|=\|Pz_{1}+Pz_{2}\|=\|Px\|>0$, and
hence $Pz_{1}$ and $Pz_{2}$ cannot both be zero. We assume $Pz_{1}\neq0$,
and then, by strict convexity of $\mathcal{H}$, $Pz_{2}=\lambda Pz_{1}$
for some $\lambda\geq0$. If $\langle z_{1}|v\rangle>\|Pz_{1}\|$
or $\langle z_{2}|v\rangle>\|Pz_{2}\|$, then $\langle x|v\rangle=\langle z_{1}|v\rangle+\langle z_{2}|v\rangle>\|Pz_{1}\|+\|Pz_{2}\|=\|Px\|$,
contradicting the assumption that $\langle x|v\rangle=\|Px\|$. Hence,
since $z_{1},z_{2}\in\mathcal{L}_{v}$, we must have $\langle z_{1}|v\rangle=\|Pz_{1}\|$
and $\langle z_{2}|v\rangle=\|Pz_{2}\|$, and therefore $\langle z_{2}|v\rangle=\|Pz_{2}\|=\lambda\|Pz_{1}\|=\lambda\langle z_{1}|v\rangle$.
Now, since $\langle z_{2}|v\rangle=\lambda\langle z_{1}|v\rangle$
and $Pz_{2}=\lambda Pz_{1}$, we obtain $z_{2}=\langle z_{2}|v\rangle v+Pz_{2}=\lambda z_{1}$,
and hence $x=z_{1}+z_{2}=(1+\lambda)z_{1}$. We conclude that $z_{1},z_{2}\in\{\lambda x:\lambda\geq0\}$.\end{proofof}
\begin{lem}
Let $\mathcal{H}$ be a real Hilbert space endowed with a Lorentz
cone $\mathcal{L}_{v}$ where $v\in\mathcal{H}$ is such that $\|v\|=1$.
If $x\in\mathcal{L}_{v}$ is such that $\langle x|v\rangle=\|Px\|$
and $0\leq y\leq x$, then $y\in\{\lambda x:\lambda\in[0,1]\}$.\end{lem}
\begin{proofof}
Since $\mathcal{L}_{v}$ is proper, this is clear if $x=0$. If $0\leq y\leq x\neq0$,
then $x=y+(x-y)$ with $y,(x-y)\in\mathcal{L}_{v}$, so that by the
previous lemma $y=\lambda x$ for some $\lambda\geq0$. If $\lambda>1$,
then $x\leq\lambda x=y\leq x$, since $\mathcal{L}_{v}$ is proper
and $y=\lambda x$, implies $x=y=0$ contradicting the assumption
$x\neq0$. We conclude that $\lambda\in[0,1]$. \end{proofof}
\begin{prop}
\label{prop:Hilbert-sp-with-lorenz-not-banach-lattice}Let $\mathcal{H}$
be a real Hilbert space endowed with a Lorentz cone $\mathcal{L}_{v}$
where $v\in\mathcal{H}$ such that $\|v\|=1$. If $\dim(\mathcal{H})\geq3$,
then $\mathcal{H}$ is not a Riesz space \textup{(}and hence not
a Banach lattice\textup{)}.\end{prop}
\begin{proofof}
Let $P$ be the orthogonal projection onto $\{v\}^{\bot}$ and $\{v,e_{1},e_{2}\}\subseteq\mathcal{H}$
be any orthonormal set. For $t\in\mathbb{R}$, we have 
\[
\{0,2e_{1}\}\leq e_{1}+te_{2}+\sqrt{t^{2}+1}v=:z_{t}.
\]
We claim that each $z_{t}$ is a minimal upper bound of $\{0,2e_{1}\}$.
We have $\langle z_{t}|v\rangle=\|Pz_{t}\|$, and hence by the previous
lemma, if $\{0,2e_{1}\}\leq y\leq z_{t}$, we must have $y=\lambda z_{t}$
for some $\lambda\in[0,1]$. If $\lambda<1$, then $\lambda\sqrt{t^{2}+1}=\langle y-2e_{1}|v\rangle$
and 
\begin{eqnarray*}
\|P(y-2e_{1})\|^{2} & = & \|P(\lambda e_{1}+\lambda te_{2}+\lambda\sqrt{t^{2}+1}v-2e_{1})\|^{2}\\
 & = & \|(\lambda-2)e_{1}+\lambda te_{2}\|^{2}\\
 & = & (\lambda-2)^{2}+\lambda^{2}t^{2}\\
 & > & 1+\lambda^{2}t^{2}\\
 & > & \lambda^{2}+\lambda^{2}t^{2}.
\end{eqnarray*}
Hence $\langle y-2e_{1}|v\rangle=\lambda\sqrt{t^{2}+1}<\|P(y-2e_{1})\|$,
contradicting $2e_{1}\leq y.$ Therefore we must have $\lambda=1$,
and $y=z_{t}$, and hence $z_{t}$ is a minimal upper bound of $\{0,2e_{1}\}$
for every $t\in\mathbb{R}$. Clearly all $z_{t}$ are distinct, and
therefore there exists no supremum of $\{0,2e_{1}\}$.
\end{proofof}
Since every Hilbert space is strictly convex, and knowing that Lorentz
cones are closed proper and generating, we conclude from Theorem \ref{thm:reflexive-u-quasi-lattice}
that every Hilbert space endowed with a Lorentz cone is an $\upsilon$-quasi-lattice.
We will now proceed to show that these spaces are absolutely monotone.
Once this has been established, Theorem \ref{thm:monotone-u-quasi-implies-mu-quasi}
will imply that they are in fact $\mu$-quasi-lattices.

The following lemma will be applied in Propositions \ref{prop:Hilbert-with-lorentz-abs-monotone}
and \ref{prop:Hilbert-Lorentz-spaces-satisfy-||x||=00003D|||x|||},
which together will show that Hilbert spaces endowed with a Lorentz
cones are in fact $1$-absolutely Davies-Ng regular. 
\begin{lem}
\label{lem:projections-preserve-upperbounds}Let $\mathcal{H}$ be
a real Hilbert space endowed with a Lorentz cone $\mathcal{L}_{v}$
where $v\in\mathcal{H}$ is such that $\|v\|=1$. Let $x\in\mathcal{H}$
and $Q$ be the orthogonal projection onto $\textup{span}\{x,v\}$.
If $\{-x,x\}\leq y$, then $\{-x,x\}\leq Qy$.\end{lem}
\begin{proofof}
Let $P$ be the orthogonal projection onto $\{v\}^{\bot}$ and $Q$
the orthogonal projection onto $\textup{span}\{x,v\}$. Let $Q^{\bot}:=\mbox{id}-Q$.
We note that $\textup{ran}(\textup{Id}-P)=\textup{span}\{v\}\subseteq\textup{ran}(Q)$,
so that $\textup{Id}-P$ and $Q$ commute, and hence $P$ and $Q$
also commute. Therefore, from
\begin{eqnarray*}
\langle v|Qy\pm x\rangle & = & \langle v|Qy+Q^{\bot}y-Q^{\bot}y\pm x\rangle\\
 & = & \langle v|y\pm x\rangle-\langle v|Q^{\bot}y\rangle\\
 & = & \langle v|y\pm x\rangle\\
 & \geq & \|P(y\pm x)\|\\
 & \geq & \|QP(y\pm x)\|\\
 & = & \|P(Qy\pm Qx)\|\\
 & = & \|P(Qy\pm x)\|,
\end{eqnarray*}
we conclude that $Qy\geq\{-x,x\}$.
\end{proofof}
The following proposition, together with Theorem \ref{thm:monotone-u-quasi-implies-mu-quasi},
will show that every Hilbert space endowed with a Lorentz cone is
in fact a $\mu$-quasi-lattice.
\begin{prop}
\label{prop:Hilbert-with-lorentz-abs-monotone}A real Hilbert space
endowed with a Lorentz cone is absolutely monotone.\end{prop}
\begin{proofof}
Let $\mathcal{H}$ be a real Hilbert space ordered by a Lorentz cone
$\mathcal{L}_{v}$, where $v\in\mathcal{H}$ is a norm-one element.
Let $P$ be the orthogonal projection onto $\{v\}^{\bot}$. Let $\{-x,x\}\leq y$
and let $Q$ denote the orthogonal projection onto $V:=\textup{span}\{x,v\}$.
By Lemma \ref{lem:projections-preserve-upperbounds}, $\{-x,x\}\leq Qy$. 

If $x\in\textup{span}\{v\}$, then $Px=0$. Also $V=\textup{span}\{v\}$,
so that $PQ=0$. Therefore $\{-x,x\}\leq Qy$ implies $\langle v|Qy\pm x\rangle\geq\|PQ(y\pm x)\|=\|Px\|=0$,
and hence $\langle v|Qy\rangle\geq|\langle v|x\rangle|$. Then $\|Qy\|=|\langle v|Qy\rangle|\geq\langle v|Qy\rangle\geq|\langle v|x\rangle|=\|x\|$,
and hence $\|x\|\leq\|Qy\|\leq\|y\|$ as was to be shown.

If $x\notin\textup{span}\{v\}$, since $0\neq Px=x-\langle v|x\rangle v\in V$,
we see that 
\[
	e_{\pm}:=(\sqrt{2}\|Px\|)^{-1}(\pm Px+\|Px\|v)
\]
 are orthonormal elements of $V\cap\mathcal{L}_{v}$. We claim that
 $V\cap\mathcal{L}_{v}=\{\lambda e_{+}+\lambda'e_{-}:\lambda,\lambda'\geq0\}$.
Let $a\in V\cap\mathcal{L}_{v}$. Since $x\notin\textup{span}\{v\}$,
$0\neq Px\in V$ is orthogonal to $v$, and hence $\{Px,v\}$ is a
basis of $V$. Then, by writing $a=\alpha Px+\beta v$ for some $\alpha,\beta\in\mathbb{R}$,
we obtain $\beta=\langle a|v\rangle\geq\|Pa\|=|\alpha|\|Px\|$. Hence,
by 
\begin{eqnarray*}
\langle a|e_{\pm}\rangle & = & (\sqrt{2}\|Px\|)^{-1}\langle\alpha Px+\beta v|\pm Px+\|Px\|v\rangle\\
 & = & (\sqrt{2}\|Px\|)^{-1}(\pm\alpha\langle Px|Px\rangle+\beta\|Px\|)\\
 & = & (\sqrt{2}\|Px\|)^{-1}(\pm\alpha\|Px\|^{2}+\beta\|Px\|)\\
 & \geq & (\sqrt{2}\|Px\|)^{-1}(\pm\alpha\|Px\|^{2}+|\alpha|\|Px\|^{2})\\
 & \geq & 0,
\end{eqnarray*}
we conclude that $V\cap\mathcal{L}_{v}=\{\lambda e_{+}+\lambda'e_{-}:\lambda,\lambda'\geq0\}$.
Now $Qy\pm x\in V\cap\mathcal{L}_{v}$ implies $\langle Qy\pm x|e_{\pm}\rangle\geq0$,
so that $\langle Qy|e_{\pm}\rangle\geq|\langle x|e_{\pm}\rangle|$,
and hence $\|x\|\leq\|Qy\|\leq\|y\|$ as was to be shown.\end{proofof}
\begin{rem}
\label{rem:vertical-cut-of-lorentz-cone-banach-lattice}If $x\notin\mbox{span}\{v\}$
we note that $(V,V\cap\mathcal{L}_{v})$ in the previous proposition
is isometrically order isomorphic to the Banach lattice $\{\mathbb{R}^{2},\|\cdot\|_{2}\}$
with the standard cone through mapping $e_{+}\in V$ and $e_{-}\in V$
to $(1,0)=:e_{1}\in\mathbb{R}^{2}$ and $(0,1)=:e_{2}\in\mathbb{R}^{2}$
respectively.
\end{rem}
We can now show that real Hilbert spaces endowed with Lorentz cones
satisfy the hypotheses of Lemma \ref{lem:sufficient-conditions-for-quasilattice-positively-attained-operator-norm}:
\begin{prop}
\label{prop:Hilbert-Lorentz-spaces-satisfy-||x||=00003D|||x|||}Let
$\mathcal{H}$ be a real Hilbert space endowed with a Lorentz cone.
Then $\|x\|=\|\left\lceil x\right\rceil \|$ for all $x\in\mathcal{H}$.
Hence  $\mathcal{H}$ is $1$-absolutely conormal. \end{prop}
\begin{proofof}
Let $v\in\mathcal{H}$ be a norm one element and order $\mathcal{H}$
with the Lorentz cone $\mathcal{L}_{v}$. We again denote the projection
onto $\{v\}^{\bot}$ by $P$. Let $x\in\mathcal{H}$ be arbitrary. 

If $x\geq0$ or $x\leq0$, then, by Theorem \ref{thm:formulas1} (6)
and (8), $\left\lceil x\right\rceil =x$ or $\left\lceil x\right\rceil =-x$,
respectively, so that $\|x\|=\|\left\lceil x\right\rceil \|$. 

It remains to show that $\|x\|=\|\left\lceil x\right\rceil \|$ when
neither $x\geq0$ nor $x\leq0$. Then $x\notin\mbox{span}\{v\}$.
We define the two dimensional subspace $V:=\mbox{span}\{x,v\}$, denote
the orthogonal projection onto $V$ by $Q$, and define $Q^{\bot}:=\mbox{Id}-Q$.
By Lemma \ref{lem:projections-preserve-upperbounds}, if $\{-x,x\}\leq w$,
then $\{-x,x\}\leq Qw$.

When $w\notin V$, we see that $Q^{\bot}w\neq0$ implies 
\begin{eqnarray*}
 &  & \|w-x\|+\|w+x\|\\
 & = & \sqrt{\|Q(w-x)\|^{2}+\|Q^{\bot}(w-x)\|^{2}}+\sqrt{\|Q(w+x)\|^{2}+\|Q^{\bot}(w+x)\|^{2}}\\
 & = & \sqrt{\|Qw-x\|^{2}+\|Q^{\bot}w\|^{2}}+\sqrt{\|Qw+x\|^{2}+\|Q^{\bot}w\|^{2}}\\
 & > & \|Qw-x\|+\|Qw+x\|.
\end{eqnarray*}
We conclude that $\left\lceil x\right\rceil $ must be an element
of $V$. Furthermore, by Proposition \ref{prop:Hilbert-with-lorentz-abs-monotone}
and Theorem \ref{thm:monotone-u-quasi-implies-mu-quasi}, $\mathcal{H}$
is a $\mu$-quasi-lattice, and hence $\left\lceil x\right\rceil \in V$
is a minimal upper bound of $\{-x,x\}$.

Finally, $V$ endowed with the cone $\mathcal{L}_{v}\cap V$ is seen
to be isometrically order isomorphic to $\{\mathbb{R}^{2},\|\cdot\|_{2}\}$
with the standard cone (cf.\,Remark \ref{rem:vertical-cut-of-lorentz-cone-banach-lattice}).
Viewing $V$ as a Banach lattice, we notice that the Banach lattice
absolute value $|x|$ in $V$ is the only minimal upper bound for
$\{-x,x\}$ in $\mathcal{H}$ that is also an element of $V$. By
the argument in the previous paragraph, we conclude that $\left\lceil x\right\rceil =|x|$,
and hence that $\|\left\lceil x\right\rceil \|=\||x|\|=\|x\|$, by
invoking the Banach lattice property $\||x|\|=\|x\|$ in $V$.
\end{proofof}
We collect the results established in this section and some of their
consequences in the following theorem. We note that (8) below resolves
the question alluded to in the introduction of the existence pre-ordered
Banach spaces $X$ and $Y$, which are not Banach lattices, while
$B(X,Y)$ is absolutely monotone.
\begin{thm}
\label{thm:hilbert-lorentz-collected-properties}Let $\mathcal{H}$
be a real Hilbert space endowed with a Lorentz cone. Then:
\begin{enumerate}
\item If $\dim(\mathcal{H})\geq3$, then $\mathcal{H}$ is not a Riesz space
\textup{(}and hence not a Banach lattice\textup{)}.
\item $\mathcal{H}$ is an $\upsilon$-quasi-lattice.
\item $\mathcal{H}$ is absolutely monotone. 
\item $\mathcal{H}$ is a $\mu$-quasi-lattice \textup{(}and its $\upsilon$--
and $\mu$-quasi-lattice structures coincide\textup{)}.
\item For every $x\in\mathcal{H}$, $\|x\|=\|\left\lceil x\right\rceil \|$.
Hence $\mathcal{H}$ is $1$-absolutely conormal. 
\item $\mathcal{H}$ is $1$-absolutely Davies-Ng regular.
\item If $X$ and $Y$ are pre-ordered Banach spaces with closed cones,
with $X$ approximately $1$-absolutely conormal and $Y$ absolutely
monotone, then the operator norms of both $B(X,\mathcal{H})$ and
$B(\mathcal{H},Y)$ are positively attained, i.e., $\|T\|=\sup\{\|Tx\|:x\geq0,\ \|x\|=1\}$
for $T\in B(X,\mathcal{H})_{+}$ or $T\in B(\mathcal{H},Y)_{+}$.
In particular, if $\mathcal{H}_{1}$ is another real Hilbert space
endowed with a Lorentz cone, then the operator norm of $B(\mathcal{H},\mathcal{H}_{1})$
is positively attained.
\item If $\alpha>0$ and $X$ and $Y$ are pre-ordered Banach spaces with
closed cones, with $X$ approximately $\alpha$-absolutely conormal
and $Y$ $\alpha$-absolutely normal, then both $B(X,\mathcal{H})$
and $B(\mathcal{H},Y)$ are $\alpha$-absolutely normal. In particular,
if $\mathcal{H}_{1}$ is another real Hilbert space endowed with a
Lorentz cone, then $B(\mathcal{H},\mathcal{H}_{1})$ is absolutely
monotone.
\end{enumerate}
\end{thm}
\begin{proofof}
The assertion (1) is Proposition \ref{prop:Hilbert-sp-with-lorenz-not-banach-lattice}.
The assertion (2) follows from Theorem \ref{thm:reflexive-u-quasi-lattice}.
The assertion (3) was established in Proposition \ref{prop:Hilbert-with-lorentz-abs-monotone}
and hence (4) follows from Corollary \ref{cor:monotone-reflexive-is-a-mu-quasi-lattice}.
The assertion (5) was established in Proposition \ref{prop:Hilbert-Lorentz-spaces-satisfy-||x||=00003D|||x|||},
and hence (6) follows from Lemma \ref{lem:sufficient-conditions-for-quasilattice-positively-attained-operator-norm}.
The assertion (7) follows from (6) and  part (2) of Theorem \ref{thm:sufficient-conditions-for-positively-attained-operator-norm}.
The assertion (8) is then immediate from (6) and part (3) of Theorem
\ref{thm-(co)normality-conditions-onXandY-implies-normality-ir-B(X,Y)}.
\end{proofof}

\section*{Acknowledgements}

The author would like to thank Marcel de Jeu, Ben de Pagter, Onno
van Gaans, Anthony Wickstead and Marten Wortel for many fruitful discussions
on this topic. Thanks should also be given to the MathOverflow community,
in particular to Anton Petrunin for ideas contained in Theorem \ref{thm:reflexive-u-quasi-lattice}.

The author's research was supported by a Vrije Competitie grant of
the Netherlands Organisation for Scientific Research (NWO).

		\bibliographystyle{amsplain}
		\bibliography{bibliography}

\end{document}